\pdfoutput=1
\documentclass[11pt, oneside]{article}   	
\usepackage{geometry}
\usepackage{graphicx}			
\usepackage{amsmath,bm}				
\usepackage{amssymb}
\usepackage{xcolor}
\usepackage{rotating}
\usepackage{caption}
\usepackage{subcaption}

\usepackage{framed}

\begin{document}
\baselineskip=2pc

\vspace{.5in}

\begin{center}

{\large\bf
Third order WENO scheme on sparse grids for hyperbolic equations
\footnote{Research supported by NSF grant DMS-1620108.}
}

\end{center}

\vspace{.1in}

\centerline{Dong Lu\footnote{Department of Applied and Computational Mathematics and Statistics,
University of Notre Dame, Notre Dame, IN 46556, USA. E-mail: dlv1@nd.edu
},
Shanqin Chen\footnote{ Department of Mathematical Sciences, Indiana University South Bend, South Bend, IN 46634, USA.
E-mail: chen39@iusb.edu},
Yong-Tao Zhang\footnote{Corresponding author. Department of Applied and Computational Mathematics and Statistics,
University of Notre Dame, Notre Dame, IN 46556, USA. E-mail: yzhang10@nd.edu
}
}

\vspace{.6in}

\abstract{The weighted essentially non-oscillatory (WENO) schemes
are a popular class of high order accurate numerical methods for solving hyperbolic
partial differential equations (PDEs). The computational cost of such schemes
increases significantly when the spatial dimensions of the PDEs are high, due to large number of spatial grid points and nonlinearity of
high order accuracy WENO schemes.
How to achieve fast computations by WENO methods for high spatial dimension PDEs is a challenging and important question.
Recently, sparse-grid has become a major approximation tool for high dimensional problems.
The open question is how to design WENO computations on sparse grids such that comparable high order accuracy of WENO schemes in smooth regions and
essentially non-oscillatory stability in non-smooth regions of the solutions can still be
achieved as that for computations on regular single grids?
In this paper, we combine the third order finite difference WENO method with sparse-grid
combination technique and solve high spatial dimension hyperbolic equations on sparse grids.
WENO interpolation is proposed for the prolongation part in sparse grid
combination techniques to deal with discontinuous solutions of hyperbolic equations.
Numerical examples are presented to show that significant computational times are saved while
both high order accuracy and stability of the WENO scheme are maintained for simulations on sparse grids.}

\vfill


\pagenumbering{arabic}

\newpage

\section{Introduction}

High order accuracy numerical methods are especially efficient for solving partial differential equations (PDEs) which contain complex solution structures.
High order numerical schemes have been applied extensively in computational fluid dynamics for solving convection dominated problems with both discontinuities / sharp gradient regions and complicated smooth structures, for example, the Rayleigh-Taylor instability simulations \cite{Remacle,SZS,ZSSZ,ZSZ}, the shock vortex interactions \cite{GP,ZZS1,ZZS2,ZZS3}, and direct simulation of compressible turbulence \cite{TWM}. Its resolution power over the lower order schemes was verified in these applications.
For hyperbolic PDEs or convection dominated problems, their solutions can develop singularities such as discontinuities, sharp gradients, discontinuous derivatives etc.
For problems containing both singularities and complicated smooth solution structures, schemes with uniform high order of accuracy in smooth regions of the solution which can also resolve singularities in an accurate and essentially non-oscillatory (ENO) fashion are desirable, since a straightforward high order approximation for the non-smooth region of a solution will generate instability called Gibbs phenomena. A popular class of such schemes is the class of weighted essentially non-oscillatory (WENO) schemes.

WENO schemes are designed based on the successful
ENO schemes \cite{Harten, SO, SO2} with additional advantages.
The first WENO scheme was constructed
by Liu, Osher, and Chan in their pioneering paper \cite{LOC} for a third order
finite volume version. In \cite{JS}, Jiang and Shu constructed
arbitrary order accurate finite
difference WENO schemes for efficiently computing multidimensional problems, with a general
framework for the design of the smoothness indicators and nonlinear
weights. To deal with complex domain geometries, WENO schemes on unstructured meshes were developed, in e.g. \cite{HS,ZS,LNSZ,DK2,ZS2,LZ}.
The main idea of the WENO schemes is to form a weighted combination
of several local reconstructions based on different stencils (usually
referred to as small stencils) and use it as the final WENO
reconstruction.  The combination coefficients (also called
nonlinear weights) depend on the linear weights, often chosen
to increase the order of accuracy over that on each small stencil, and on the
smoothness indicators which measure the smoothness of the
reconstructed function in the relevant small stencils.
Hence an adaptive approximation or reconstruction procedure is
actually the essential part of the WENO schemes.

Since WENO schemes were designed to deal with problems with complicated solution structures,
their sophisticated nonlinear properties and high order accuracy requires more operations than
many other schemes. For PDEs with high spatial dimensions, large number of spatial grid points
leads to significant increase of the computational cost for WENO schemes, especially for long time simulations or steady state computations \cite{ZZQ,ZZC,XZZS,HLT2,HHSSXZ,WZZS}. It is challenging and important
to achieve fast computations when WENO methods are used for solving high spatial dimension PDEs.

In recent years, sparse-grid techniques have been used broadly as an efficient approximation tool for high-dimensional problems
in many scientific and engineering applications. Discretizations on sparse grids involve
$O(N\cdot (\log N)^{d-1})$ degrees of freedom only, where $d$ denotes the dimensionality of the underling problems  and
$N$ is the number of grid points in one coordinate direction. A detailed review on sparse-grid technique can be found in
\cite{BG}. Sparse-grid techniques were introduced by Zenger \cite{Zg} in 1991 to reduce the number of degrees of freedom in finite element calculations. The sparse-grid combination technique, which was introduced in 1992 by Griebel et al. \cite{GSZ}, can be seen as a practical implementation of the sparse-grid technique. In the sparse-grid combination technique, the final solution is a linear combination of solutions on semi-coarsened grids, where the coefficients of the combination are chosen such that there is a canceling in leading-order error terms and the accuracy order can be kept to be the same as that on single full grids \cite{LKV1,LKV2,GSZ}.
Recently in \cite{LuZhang1}, the sparse-grid combination technique has been used in Krylov implicit integration factor methods \cite{CZ,JiangZhang,JiangZhang2,LuZhang2}
to efficiently solve high spatial dimension convection-diffusion  equations.

Our goal is to apply sparse-grid techniques in high order WENO schemes to achieve more efficient computations than that in their regular performance in solving multidimensional PDEs. The open question is how to design WENO computations on sparse grids such that comparable high order accuracy of WENO schemes in smooth regions and
essentially non-oscillatory stability in non-smooth regions of the solutions can still be
achieved as that for computations on regular single grids? This is not straightforward due to the high nonlinearity of high order WENO schemes.
In this paper, we design and test a third order sparse grid WENO finite difference scheme for solving hyperbolic PDEs by using the sparse-grid combination approach.
To deal with discontinuous solutions of hyperbolic PDEs, we apply WENO interpolation for the prolongation part in sparse-grid
combination techniques. Both two dimensional (2D) and three dimensional (3D) numerical examples with smooth or non-smooth solutions are presented to show that significant computational times are saved, while both accuracy and stability of the WENO scheme are maintained for simulations on sparse grids.
The rest of the paper is organized as following. In Section 2, we design the third order sparse grid WENO scheme. In Section 3, numerical experiments
are presented to test the sparse grid WENO method and show significant savings in
computational costs by comparisons with single-grid computations. Conclusions are given in Section 4.

\section{A third order sparse grid WENO finite difference scheme}
We consider multidimensional hyperbolic PDEs
\begin{equation}
\label{eq1}
{u}_t+\nabla\cdot\vec{f}({u})=0,
\end{equation}
where ${u}(\vec x,t)$ is the unknown, and $\vec{f}=({f}_1, \cdots, {f}_d)^T$ is the vector of
flux functions in $d$ spatial dimensions respectively.
%For simplicity of presentation, we first describe the algorithm in details using 2D problems. Algorithm for higher dimensional problems is similar and formulas will be %given. Numerical experiments are performed for both 2D and 3D problems.
The method of lines (MOL)
 is applied to the equation (\ref{eq1}). The third order finite difference WENO scheme
with Lax-Friedrichs flux splitting is used for spatial discretizations. In this section,
we first describe the spatial discretization by finite difference WENO scheme, then
the sparse-grid combination approach with the WENO prolongation is introduced and a complete algorithm is given.

\subsection{WENO discretization}
For the hyperbolic PDEs (\ref{eq1}), the conservative finite difference scheme we use approximates the point values at a
uniform (or smoothly varying) grid in a conservative fashion. Since the finite difference WENO schemes
approximate derivatives of multi-dimension in a dimension by dimension way, we will just describe the discretization of derivatives for
one spatial direction. As a general notation,
for example, we consider the $x$-direction derivative $f(u)_x$. Its value at a grid point with $x$-coordinate $x_i$ on a uniform grid with $x$-direction grid size $\Delta x$
is approximated by a conservative flux difference
\begin{equation}
\label{eq2.9}
f(u)_x|_{x=x_i}\approx \frac{1}{\Delta x}(\hat f_{i+1/2}-\hat f_{i-1/2}),
\end{equation}
where for the third order WENO scheme the numerical flux $\hat
f_{i+1/2}$ depends on the three-point values $f(u_{l})$ (here for
the simplicity of notations, we use $u_l$ to denote the value of
the numerical solution $u$ at the point $x=x_l$ along the lines of other spatial directions, e.g.,
$y=y_j, z=z_k,$ etc, with the understanding that the value could be
different for different coordinates of other spatial directions), $l=i-1,i,i+1$,
when the wind is positive (i.e., when $f'(u)\geq 0$ for the scalar
case, or when the corresponding eigenvalue is positive for the
system case with a local characteristic decomposition). This
numerical flux $\hat f_{i+1/2}$ is written as a convex combination
of two second order numerical fluxes based on two different
substencils of two points each, and the combination coefficients
depend on a ``smoothness indicator'' measuring the smoothness of
the solution in each substencil. The detailed formula is
\begin{equation}
\label{eq2.7}
\hat f_{i+1/2}=w_0\left[\frac{1}{2}f(u_i)+\frac{1}{2}f(u_{i+1})\right]
+w_1\left[-\frac{1}{2}f(u_{i-1})+\frac{3}{2}f(u_{i})\right],
\end{equation}
where
\begin{equation}
w_r=\frac{\alpha_r}{\alpha_1+\alpha_2}, \qquad
\alpha_r=\frac{d_r}{(\epsilon+\beta_r)^2}, \qquad r=0, 1.
\label{eq2.6}
\end{equation}
$d_0=2/3, d_1=1/3$ are called the ``linear weights'', and
$\beta_0=(f(u_{i+1})-f(u_i))^2, \beta_1=(f(u_i)-f(u_{i-1}))^2$ are
called the ``smoothness indicators''. $\epsilon$ is a small
positive number chosen to avoid the denominator becoming $0$.

When the wind is negative (i.e., when $f'(u)< 0$), right-biased
stencil with numerical values $f(u_i), f(u_{i+1})$ and
$f(u_{i+2})$ are used to construct a third order WENO
approximation to the numerical flux $\hat f_{i+1/2}$. The formulae
 for negative and positive wind cases are symmetric with respect to the
point $x_{i+1/2}$. For the general case of $f(u)$, we perform the
``Lax-Friedrichs flux splitting''
\begin{equation}
f^+(u)=\frac{1}{2}(f(u)+\alpha u), \qquad
f^-(u)=\frac{1}{2}(f(u)-\alpha u),
\end{equation}
where $\alpha=\max_u|f'(u)|$. $f^+(u)$ is the positive wind part,
and $f^-(u)$ is the negative wind part. Corresponding WENO
approximations are applied to find numerical fluxes $\hat
f^+_{i+1/2}$ and $\hat f^-_{i+1/2}$ respectively. Similar
procedures are applied to other spatial directions' derivatives. See \cite{Shu,ZSNA} for more details.

\subsection{WENO scheme on sparse grids}

To efficiently solve the hyperbolic equations (\ref{eq1}) on high spatial dimensions by WENO schemes, we study the WENO schemes on sparse grids by sparse-grid combination technique. In this paper, we focus on the third order finite difference WENO scheme given in the last section.

The basic idea of sparse-grid combination technique is that by combining several solutions on different semi-coarsened grids (sparse grids), a final solution on the most refined mesh is obtained. The most refined mesh is corresponding to the usual single full grid. Since the PDEs are solved on semi-coarsened grids which have much fewer grid points than the single full grid, computation costs are saved a lot. The final solution obtained by sparse-grid combination technique is required to
have comparable accuracy to that by solving the PDEs directly on a single full grid. For example see \cite{GSZ,LKV1,LKV2,LuZhang1}.

\begin{figure}
\centering
\includegraphics[width = 4in]{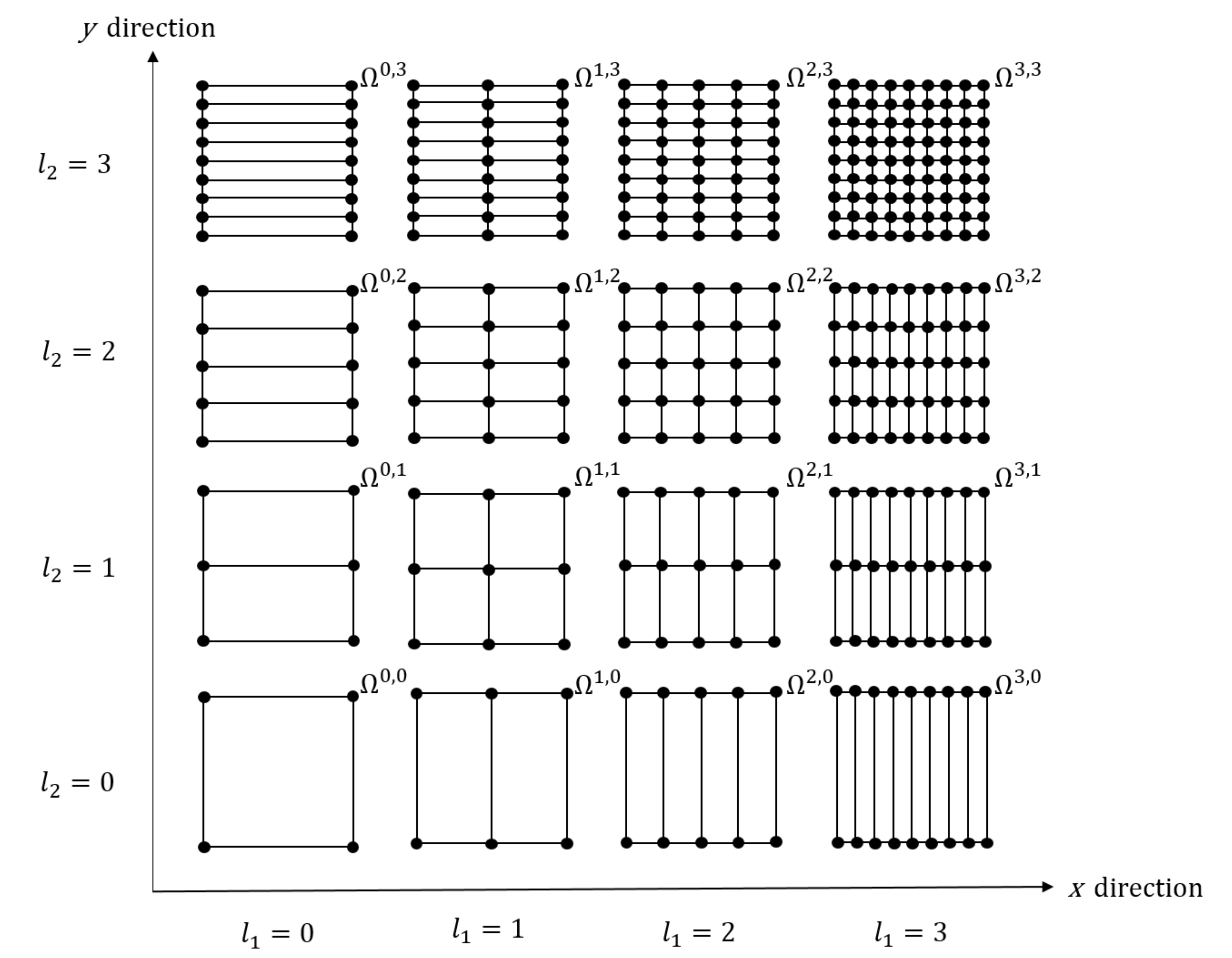}
\caption{Semi-coarsened sparse grids $\{\Omega^{l_1,l_2}\}$ with the finest level $N_L=3$.}
\label{sparse_grids}
\end{figure}

2D case is used here as the example to illustrate the idea. Algorithms are similar for higher dimensional cases. Let's consider a 2D domain $[a, b]^2$.
The semi-coarsened grids are constructed as follows. First the domain is partitioned into the coarsest mesh, which is called a root grid $\Omega^{0,0}$ with $N_r$ cells in each direction. The root grid mesh size is $H=\frac{b-a}{N_r}$. The multi-level refinement on the root grid is performed to obtain a family of semi-coarsened grids \{$\Omega^{l_1,l_2}$\}. The semi-coarsened grid \{$\Omega^{l_1,l_2}$\} has mesh sizes $h_{l_1}=2^{-l_1}H$ in the $x$ direction and $h_{l_2}=2^{-l_2}H$ in the $y$ direction, where $l_1=0,1,\cdots,N_L$, $l_2=0,1,\cdots,N_L$ (see Figure \ref{sparse_grids}). Superscripts $l_1, l_2$ denote the level of refinement relative to the root grid $\Omega^{0,0}$, and $N_L$ denotes the finest level. Therefore, our finest grid is $\Omega^{N_L,N_L}$ with mesh size $h=2^{-N_L}H$ for both $x$ and $y$ directions.

To solve equation (\ref{eq1}), we use the third order WENO scheme described in Section 2.1 for spatial discretizations with the third order TVD Runge-Kutta scheme \cite{SO,Shu} for time discretization.  Following the spare-grid combination techniques, rather than on a single full grid, the PDE (\ref{eq1}) is solved on the following $(2N_L+1)$ sparse grids $\{\Omega^{l_1,l_2}\}_I$:
\[
\Big\{
\Omega^{0,N_L}, \Omega^{1,N_L-1}, \cdots, \Omega^{N_L-1,1}, \Omega^{N_L,0}
\Big\}\quad
\text{and} \quad
\Big\{
\Omega^{0,N_L-1}, \Omega^{1,N_L-2}, \cdots, \Omega^{N_L-2,1}, \Omega^{N_L-1,0}
\Big\}.
\]
And $I$ denotes the index set
\[
I=\big\{(l_1,l_2)| l_1+l_2=N_L \quad \text{or} \quad l_1+l_2=N_L-1\big\}.
\]
By carrying out time marching of the PDE using Runge-Kutta WENO scheme on these $(2N_L+1)$ sparse grids, we obtain $(2N_L+1)$ sets of numerical solutions $\{U^{l_1,l_2}\}_I$, where one set of numerical solutions is obtained on each sparse grid. The next step is to combine solutions on sparse grids to obtain the final solution on the finest grid $\Omega^{N_L,N_L}$. The key point here is that the PDE is never solved directly on $\Omega^{N_L,N_L}$ in order to save computational costs. In order to obtain numerical solutions on the finest grid $\Omega^{N_L,N_L}$,
we apply a prolongation operator $P^{N_L,N_L}$, which will be discussed in the following subsections, on each sparse grid solution $U^{l_1,l_2}$ to obtain $(2N_L+1)$ solutions on the finest grid. And finally, these solutions are combined to form the final solution $\hat{U}^{N_L,N_L}$ on $\Omega^{N_L,N_L}$.

\subsubsection{Lagrange prolongation}

We provide details on the prolongation operator $P^{N_L,N_L}$.
Prolongation operator $P^{N_L,N_L}$ maps numerical solutions $\{U^{l_1,l_2}\}_I$ on sparse grids onto the finest grid $\Omega^{N_L,N_L}$. And a prolongation operator is basically an interpolation operator. For example, $U^{l_1,l_2}$ is numerical solution on $\Omega^{l_1,l_2}$, then $P^{N_L,N_L}U^{l_1,l_2}$ gives numerical values on the most refined mesh $\Omega^{N_L,N_L}$. For smooth solutions, the regular Lagrange interpolation can be used directly. The interpolations are performed in the dimension by dimension way.
For the 2D case, first in grid lines of one direction (e.g. the $x$ direction with a fixed $y$-coordinate), we construct $(N_r2^{l_1-1})$ quadratic interpolation polynomials $P_i^2(x)$, $i=1,\cdots,N_r2^{l_1-1}$, by the third order Lagrange interpolation. Each interpolation uses three adjacent grid points to construct a quadratic polynomial. Note that
a higher order interpolation is needed for comparable numerical accuracy as that of the numerical schemes, if higher order accuracy numerical schemes are used to solve PDEs on sparse grids (see \cite{GSZ,LKV1,LKV2}).
Then we evaluate $P_i^2(x)$ on the grid points of $\Omega^{N_L,l_2}$, which is the most refined meshes in the $x$ direction. Next, in grid lines of the other direction (e.g. the $y$ direction with a fixed $x$-coordinate), we construct $(N_r2^{l_2-1})$ quadratic interpolation polynomials $P_j^2(y)$, $j=1,\cdots,N_r2^{l_2-1}$, and evaluate them on the grid points of $\Omega^{N_L,N_L}$. Then we get $P^{N_L,N_L}U^{l_1,l_2}$, defined on the finest grid $\Omega^{N_L,N_L}$.

\subsubsection{WENO prolongation}

In general, since the solutions may develop discontinuities, instead of Lagrange interpolation, more robust WENO interpolations are used in the prolongation. Replacing
the third order Lagrange interpolation by a third order WENO interpolation in the procedure described in the section 2.2.1, we obtain a third order WENO prolongation in the sparse-grid combination technique. We provide the detailed formulas for a third order WENO interpolation here. Given numerical values $u_{i-1}$, $u_i$ and $u_{i+1}$ at the
grid points $x_{i-1}$, $x_i$ and $x_{i+1}$, we find a third order WENO interpolation $u_{WENO}(x)$ for any point $x \in [x_{i-1/2}, x_{i+1/2})$. Here $x_{i-1/2}=(x_{i-1}+x_i)/2$ and $x_{i+1/2}=(x_{i}+x_{i+1})/2$. Denote the uniform grid size by $h$, then $\forall x \in [x_{i-1/2}, x_{i+1/2})$, $x=x_{i-1}+\tilde\alpha h$ where
$\tilde\alpha \in [1/2, 3/2)$. The third order WENO interpolation $u_{WENO}(x)$ at a point $x$ is
\begin{equation}
u_{WENO}(x) = w_1 P_{(1)}^1(x) + w_2 P_{(2)}^1(x),
\label{wp1}
\end{equation}
where $P_{(1)}^1(x)$ and $P_{(2)}^1(x)$ are two second order approximations
\begin{equation}
P_{(1)}^1(x) = \tilde\alpha u_i -(\tilde\alpha - 1) u_{i-1},
\qquad
P_{(2)}^1(x) = (\tilde\alpha - 1) u_{i+1} - (\tilde\alpha - 2) u_{i}.
\label{wp2}
\end{equation}
$w_1$ and $w_2$ are nonlinear weights
\begin{equation}
w_1 = \frac{\tilde w_1}{\tilde w_1 + \tilde w_2}, \qquad
w_2 = 1 - w_1,
\label{wp3}
\end{equation}
where
\begin{equation}
\tilde w_1 = \frac{\gamma_1}{(\epsilon + \beta_1)^2}, \qquad
\tilde w_2 = \frac{\gamma_2}{(\epsilon + \beta_2)^2},
\label{wp4}
\end{equation}
and $\gamma_1=1-\tilde\alpha/2$, $\gamma_2=\tilde\alpha/2$, $\beta_1 = (u_i-u_{i-1})^2$, $\beta_2 = (u_{i+1}-u_{i})^2$.
Again, $\epsilon$ is a small
positive number chosen to avoid the denominator becoming $0$.

\subsubsection{Algorithm}

We summarize the algorithm of WENO scheme on sparse grids as following.

%\begin{framed}
\bigskip
\noindent\textbf{Algorithm: WENO scheme with sparse-grid combination technique}
\begin{itemize}
\item Step 1: Restrict the initial condition $u(x,y,0)$ to $(2N_L+1)$ sparse grids $\{\Omega^{l_1,l_2}\}_I$ defined above. Here ``Restrict'' means that functions are
evaluated at grid points;
\item Step 2: On each sparse grid $\Omega^{l_1,l_2}$, solve the equation (\ref{eq1}) by Runge-Kutta WENO scheme to reach the final time $T$. Then we get $(2N_L+1)$ sets of solutions $\{U^{l_1,l_2}\}_I$;
\item Step 3: At the final time $T$,
    \begin{itemize}
    \item on each grid $\Omega^{l_1,l_2}$, apply prolongation operator $P^{N_L,N_L}$ on $U^{l_1,l_2}$. Then we get $P^{N_L,N_L}U^{l_1,l_2}$, defined on the most refined mesh $\Omega^{N_L,N_L}$. For smooth solutions, the regular Lagrange prolongation can be used directly. In general, WENO prolongation is used;
    \item do the combination to get the final solution
    \begin{equation}
    \hat{U}^{N_L,N_L}=\sum_{l_1+l_2=N_L}P^{N_L,N_L}U^{l_1,l_2}-\sum_{l_1+l_2=N_L-1}P^{N_L,N_L}U^{l_1,l_2}.
    \end{equation}
    \end{itemize}
\end{itemize}
%\end{framed}

For three dimensional (3D) or higher dimensional problems, the algorithm is similar although prolongation operations are performed in additional spatial directions. The sparse-grid combination formula for higher dimensional cases can be found in the literature (e.g. \cite{GSZ}). Specifically the 3D formula is
\begin{equation}
\begin{aligned}
&\hat{U}^{N_L,N_L,N_L}=\sum_{l_1+l_2+l_3=N_L}P^{N_L,N_L,N_L}U^{l_1,l_2,l_3}-2\sum_{l_1+l_2+l_3=N_L-1}P^{N_L,N_L,N_L}U^{l_1,l_2,l_3}\\
&+\sum_{l_1+l_2+l_3=N_L-2}P^{N_L,N_L,N_L}U^{l_1,l_2,l_3}.\\
\end{aligned}
\end{equation}

\bigskip
\noindent{\bf Remark:} Linear error analysis of the sparse-grid combination technique for a linear advection equation solved by an upwind scheme has been performed in
\cite{LKV1}. In this paper, we focus on the algorithm development and its numerical experiments for the nonlinear WENO scheme on sparse grids. The nonlinear analysis of the
scheme will be performed in one of our future work.

%\newpage
\section{Numerical Experiments}
In this section, we use various numerical examples to show the computational efficiency of the third order WENO scheme with sparse-grid combination technique on sparse grids, by comparing to the same scheme on regular grids. Examples include both 2D and 3D numerical examples with smooth or non-smooth solutions. For each example, we compute numerical accuracy errors and convergence orders of the schemes, and record CPU times. Here in the data Tables and texts of this section, $N_h\times N_h$
denotes the most refined mesh $\Omega^{N_L,N_L}$ in sparse grids or a regular mesh in single grid computations.

The third order linear scheme is obtained by replacing nonlinear weights $w_0$ and $w_1$ in WENO approximation (\ref{eq2.7}) with linear weights $d_0$ and $d_1$.
Linear schemes are stable and efficient for solving problems with smooth solutions, and they serve as the base schemes for high order WENO schemes. We also test
the computational efficiency of the third order linear scheme on sparse grids for solving problems with smooth solutions.

For computations on sparse grids, PDEs are evolved on different semi-coarsened sparse grids. How to choose time step sizes for
each individual time evolution is an interesting question. Via numerical experiments, we find that time step sizes on all semi-coarsened sparse grids need to take the same value. It is
determined by the spatial grid size $h$ of the most refined grid $\Omega^{N_L,N_L}$ and the chosen CFL number. Numerical experiments show that the desired numerical accuracy are reached with time step sizes taken this way.
Hence for a general problem, the numerical experiments in this paper suggest that time step sizes on all semi-coarsened sparse grids should be determined by the spatial grid size $h$ of the most refined grid $\Omega^{N_L,N_L}$.
All of the numerical simulations in this paper are performed on a 2.3 GHz, 16GB RAM Linux workstation.

\subsection{Examples with Smooth Solutions}

In this section, we test the scheme on sparse grids for solving problems with smooth solutions.

\paragraph{Example 1 (A 3D Linear equation):}
\begin{equation}
\label{eqn:3.20}
\begin{cases}
u_t+u_x+u_y+u_z=0, \qquad -2 \leq x \leq 2, -2 \leq y \leq 2, -2 \leq z \leq 2;\\
u(x,y,z,0) =\sin(\frac{\pi}{2}(x+y+z)),
\end{cases}
\end{equation}
with periodic boundary condition. We compute this 3D problem till final time $T=1$ by both the third order linear scheme and WENO scheme
on both single grids and sparse grids, and compare their computational efficiency. The $L^\infty$ errors, $L^2$ errors, the corresponding numerical accuracy orders, and CPU times on
successively refined meshes to show the efficiency of computations on sparse grids are reported. To refine meshes for computations on sparse grids, we
refine the root grid $\Omega^{0,0,0}$, and
keep the number of semi-coarsened sparse-grid levels (total $N_L+1$ levels) unchanged. For example, sparse-grid with a $10 \times 10\times 10$ root grid and
$N_L=3$ has the finest mesh $80 \times 80\times 80$. If the root grid is refined once to be $20 \times 20\times 20$, with $N_L=3$ unchanged we
can obtain the finest mesh $160 \times 160\times 160$. The numerical errors, accuracy orders, and CPU times are listed in Table \ref{tab:ln_ln_3d} for the third order linear scheme and Table \ref{tab:ln_weno_3d} for the third order WENO scheme. Two different finest levels $N_L=3$ and $N_L=2$ are tested in sparse-grid computations. From Table \ref{tab:ln_ln_3d}, we can see that for the linear scheme, the computations on single grids and sparse grids achieve
the comparable numerical errors and the third order accuracy. However, computations on sparse-grid are much more efficient than those on single-grid.
Comparing the CPU times in Table\ref{tab:ln_ln_3d}, we can see that for computations on sparse grids with $N_L=3$, more than $80\%$ computation time can be saved
to reach the similar error levels as that on single grids. If $N_L=2$, $55\% \sim 64\%$ computation time is saved. From Table \ref{tab:ln_weno_3d} for the third order WENO scheme, we can see that on relatively coarse meshes, the numerical errors of computations on sparse grids are larger than that on single grids. However, with more refined meshes the computations on sparse grids show superconvergence and achieve comparable numerical errors and accuracy as that on single grids. For $N_L=3$, $84\%$ CPU time can be saved for the computation on the $640\times 640\times 640$ mesh. While $N_L=2$, $63\%$ CPU time is saved on the $320\times 320\times 320$ mesh and
$67\%$ CPU time is saved on the $640\times 640\times 640$ mesh. If we compare the results of $N_L=2$ and $N_L=3$ in Table \ref{tab:ln_weno_3d} for the WENO scheme, we find
that the computations on sparse grids can achieve smaller numerical errors with a smaller finest-level $N_L$, while CPU costs are less if $N_L$ is larger.

\begin{table}[htbp]\footnotesize
\centering
\begin{tabular}{c c c c c c c c}
\hline
\multicolumn{8}{c}{ Single-grid}\\ \hline
& & $N_h\times N_h\times N_h$	&$L^\infty$ error		&Order	&$L^2$ error &Order &CPU(s)	\\
\hline
&&$80\times 80\times 80$          &$2.30\times10^{-4}$ &       &$1.63\times 10^{-4}$    &    &52.95	    	\\
&&$160\times 160\times 160$       &$2.88\times10^{-5}$ &3.00   &$2.04\times 10^{-5}$    &3.00&930.30       	\\
&&$320\times 320\times 320$       &$3.60\times10^{-6}$ &3.00   &$2.55\times 10^{-6}$    &3.00&15,030.00    	\\
&&$640\times 640\times 640$       &$4.50\times10^{-7}$ &3.00   &$3.18\times 10^{-7}$    &3.00     &261,972.90 		\\
\hline
\multicolumn{8}{c}{ Sparse-grid, refine root grids, $N_L=3$}\\ \hline
$N_r$   &$N_L$  &$N_h\times N_h\times N_h$	&$L^\infty$ error   &Order &$L^2$ error   &Order	&CPU(s)	      \\
\hline
10  &3  &$80\times 80\times 80$          &$7.16\times10^{-4}$ &       &$4.91\times 10^{-4}$    &    &9.87	    	\\
20  &3  &$160\times 160\times 160$       &$3.96\times10^{-5}$ &4.17   &$2.79\times 10^{-5}$    &4.14&148.34     	\\
40  &3  &$320\times 320\times 320$       &$3.88\times10^{-6}$ &3.35   &$2.74\times 10^{-6}$    &3.35&3,087.23    	\\
80  &3  &$640\times 640\times 640$       &$4.58\times10^{-7}$ &3.08   &$3.24\times 10^{-7}$    &3.08&52,702.50		\\
\hline
\multicolumn{8}{c}{ Sparse-grid, refine root grids, $N_L=2$}\\ \hline
$N_r$   &$N_L$  &$N_h\times N_h\times N_h$	&$L^\infty$ error   &Order &$L^2$ error    &Order	&CPU(s)	       \\
\hline
20  &2  &$80\times 80\times 80$         &$2.32\times10^{-4}$ &       &$1.64\times 10^{-4}$    &     &18.97	    	\\
40  &2  &$160\times 160\times 160$       &$2.85\times10^{-5}$ &3.03   &$2.02\times 10^{-5}$    &3.02&368.89      	\\
80  &2  &$320\times 320\times 320$       &$3.59\times10^{-6}$ &2.99   &$2.54\times 10^{-6}$    &2.99&6,741.18     	\\
160  &2  &$640\times 640\times 640$      &$4.49\times10^{-7}$ &3.00   &$3.18\times 10^{-7}$    &3.00&106,619.00  	\\
\hline
\end{tabular}
\caption{\footnotesize{Example 1, Linear scheme,
comparison of numerical errors and CPU times for computations on single-grid and sparse-grid. Lagrange interpolation for prolongation is used in sparse-grid computations.
Final time $T=1.0$. CFL number is $0.75$.
$N_r$: number of cells in each spatial direction of a root grid.
$N_L$: the finest level in a sparse-grid computation.
CPU: CPU time for a complete simulation. CPU time unit: seconds.}}
\label{tab:ln_ln_3d}
\end{table}

\begin{table}[htbp]\footnotesize
\centering
\begin{tabular}{c c c c c c c c}
\hline
\multicolumn{8}{c}{ Single-grid}\\ \hline
& & $N_h\times N_h\times N_h$	&$L^\infty$ error		&Order	&$L^2$ error &Order &CPU(s)	\\
\hline
&&$80\times 80\times 80$          &$8.30\times10^{-4}$ &       &$3.81\times 10^{-4}$    &    &90.74	    	\\
&&$160\times 160\times 160$       &$4.83\times10^{-5}$ &4.10   &$2.56\times 10^{-5}$    &3.89&1,561.30      	\\
&&$320\times 320\times 320$       &$4.21\times10^{-6}$ &3.52   &$2.67\times 10^{-6}$    &3.26&30,656.40     	\\
&&$640\times 640\times 640$       &$4.69\times10^{-7}$ &3.17   &$3.22\times 10^{-7}$    &3.05&521,562.00  	\\
\hline
\multicolumn{8}{c}{ Sparse-grid, refine root grids, $N_L=3$}\\ \hline
$N_r$   &$N_L$  &$N_h\times N_h\times N_h$	&$L^\infty$ error   &Order &$L^2$ error   &Order	&CPU(s)	      \\
\hline
10  &3  &$80\times 80\times 80$          &$1.35\times10^{-1}$ &       &$5.82\times 10^{-2}$    &    &18.15	     	\\
20  &3  &$160\times 160\times 160$       &$1.46\times10^{-2}$ &3.21   &$6.86\times 10^{-3}$    &3.09&306.79      	\\
40  &3  &$320\times 320\times 320$       &$1.64\times10^{-4}$ &6.47   &$7.93\times 10^{-5}$    &6.44&5,325.82	  	\\
80  &3  &$640\times 640\times 640$       &$6.97\times10^{-7}$ &7.88   &$3.40\times 10^{-7}$    &7.86&82,575.40	 	\\
\hline
\multicolumn{8}{c}{ Sparse-grid, refine root grids, $N_L=2$}\\ \hline
$N_r$   &$N_L$  &$N_h\times N_h\times N_h$	&$L^\infty$ error   &Order &$L^2$ error    &Order	&CPU(s)	       \\
\hline
20  &2  &$80\times 80\times 80$         &$2.77\times10^{-2}$ &        &$1.38\times 10^{-2}$    &     &33.49	      	\\
40  &2  &$160\times 160\times 160$       &$6.69\times10^{-4}$ &5.37   &$3.09\times 10^{-4}$    &5.48 &777.97      	\\
80  &2  &$320\times 320\times 320$       &$5.47\times10^{-6}$ &6.93   &$2.71\times 10^{-6}$    &6.83 &11,242.90	  	\\
160  &2  &$640\times 640\times 640$      &$4.70\times10^{-7}$ &3.55   &$3.21\times 10^{-7}$    &3.08 &171,576.00  	\\
\hline
\end{tabular}
\caption{\footnotesize{Example 1, WENO scheme,
comparison of numerical errors and CPU times for computations on single-grid and sparse-grid. Lagrange interpolation for prolongation is used in sparse-grid computations.
Final time $T=1.0$. CFL number is $0.75$.
$N_r$: number of cells in each spatial direction of a root grid.
$N_L$: the finest level in a sparse-grid computation.
CPU: CPU time for a complete simulation. CPU time unit: seconds.}}
\label{tab:ln_weno_3d}
\end{table}

\paragraph{Example 2 (A 2D Nonlinear equation):}
\begin{equation}
\label{eqn:Ex2}
\begin{cases}
u_t+(\frac{1}{2}u^2)_x+(\frac{1}{2}u^2)_y=-0.1u, \qquad 0 \leq x \leq 2\pi, 0 \leq y \leq 2\pi;\\
u(x,y,0)=\sin(x-y),
\end{cases}
\end{equation}
with periodic boundary condition. The exact solution of the problem is $u(x,y,t)=e^{-0.1t}\sin(x-y)$.
We compute this 2D nonlinear problem till final time $T=1$ by both the third order linear scheme and WENO scheme
on both single grids and sparse grids, and compare their computational efficiency. In Table \ref{tab:nonln_ln_2d} and Table \ref{tab:nonln_weno_2d},
 the $L^\infty$ errors, $L^2$ errors, the corresponding numerical accuracy orders, and CPU times on
successively refined meshes to show the efficiency of computations on sparse grids are reported. Similar as Example 1, for the linear scheme, comparable numerical errors and third order accuracy are obtained on single grids and sparse grids, except that the errors are larger on the $80 \times 80$ mesh for the $N_L=3$ case.
For the WENO scheme, numerical errors and accuracy on sparse grids are comparable to those on single grids if the mesh is relatively refined. Comparing the CPU costs of computations on sparse grids and single grids, we find that for 2D problem, the saving of CPU times of sparse-grid computations is less than that for 3D problem. As that shown in Table \ref{tab:nonln_ln_2d} and Table \ref{tab:nonln_weno_2d}, about $30\%$ computation time can be saved if $N_L=3$ in sparse-grid computations. The CPU times of sparse-grid computations with $N_L=2$ are similar as that of single-grid computations.
Again similar as Example 1, the computations on sparse grids can achieve smaller numerical errors with a smaller finest-level $N_L$, while CPU costs are less if $N_L$ is larger.
In sparse-grid computations presented in this paper, the saving of CPU times of sparse-grid computations is due to the fact that less number of grid points is used.
For 2D sparse grids with $N_r$ number of cells in each spatial direction of the root grid and the finest level $N_L$, the number of
grid points of the most refined grid $\Omega^{N_L,N_L}$ is
 $(N_r \cdot 2^{N_L}+1)^2$. Note that this is also the number of
grid points which are used in discretizing and solving the PDEs on the corresponding single grids. While in
the sparse-grid computations, the number of
grid points which are used in discretizing and solving the PDEs is
$$\sum_{i=0}^{N_L}(N_r\cdot 2^i+1)\cdot(N_r \cdot 2^{N_L-i}+1)+\sum_{i=0}^{N_L-1}(N_r \cdot 2^i+1)\cdot(N_r \cdot 2^{N_L-i-1}+1).$$
Note that this count of grid points on sparse grids does not include those grid points used in the prolongation step, which is not in the time evolution process and done only at the final time step, and does not directly involve solving the PDEs.
In Table \ref{tab:gridcount_2d}, we list numbers of spatial grid points which are used in discretizing and solving this example on sparse grids and single grids.
Especially, the ratios of numbers of spatial grid points on sparse grids to that on single grids are shown. They are compared with the ratios of CPU times of the
third order WENO scheme on sparse grids to that on single grids. We obtain consistent results. If $N_L=3$, sparse-grid computations use about $70\%$ grid points of that
in single-grid computations, which leads to the saving of around $30\%$ CPU times. For $N_L=2$, the numbers of grid points used in sparse-grid and single-grid computations
are similar, hence their CPU times are also similar.

\begin{table}[htbp]\footnotesize
\centering
\begin{tabular}{c c c c c c c c c}
\hline
\multicolumn{9}{c}{ Single-grid}\\ \hline
& & $N_h\times N_h$		&$L^\infty$ error		&Order  &$L^2$ error    &Order  &CPU(s)		&Ratio	\\
\hline
&&$80\times 80$         &$6.95\times10^{-5}$ &       &$4.91\times 10^{-5}$    &    &0.40	    &	\\
&&$160\times 160$       &$8.69\times10^{-6}$ &3.00   &$6.14\times 10^{-6}$    &3.00&3.20	    &7.93	\\
&&$320\times 320$       &$1.09\times10^{-6}$ &3.00   &$7.68\times 10^{-7}$    &3.00&26.02	    &8.14	\\
&&$640\times 640$       &$1.36\times10^{-7}$ &3.00   &$9.60\times 10^{-8}$    &3.00&220.63	&8.48	\\
\hline
\multicolumn{9}{c}{ Sparse-grid, refine root grids, $N_L=3$}\\ \hline
$N_r$   &$N_L$  &$N_h\times N_h$	&$L^\infty$ error   &Order	    &$L^2$ error    &Order  &CPU(s)	     &Ratio\\
\hline
10  &3  &$80\times 80$         &$3.64\times10^{-4}$ &       &$1.57\times 10^{-4}$    &    &0.33	    &	\\
20  &3  &$160\times 160$       &$1.14\times10^{-5}$ &4.99   &$6.50\times 10^{-6}$    &4.59&2.33	    &7.16	\\
40  &3  &$320\times 320$       &$1.11\times10^{-6}$ &3.36   &$7.73\times 10^{-7}$    &3.07&17.76	&7.61 	\\
80  &3  &$640\times 640$       &$1.36\times10^{-7}$ &3.03   &$9.61\times 10^{-8}$    &3.01&141.14	&7.95	\\
\hline
\multicolumn{9}{c}{ Sparse-grid, refine root grids, $N_L=2$}\\ \hline
$N_r$   &$N_L$  &$N_h\times N_h$	&$L^\infty$ error   &Order	    &$L^2$ error    &Order  &CPU(s)	     &Ratio\\
\hline
20  &2  &$80\times 80$         &$8.05\times10^{-5}$ &       &$5.10\times 10^{-5}$    &    &0.43	    &	\\
40  &2  &$160\times 160$       &$8.76\times10^{-6}$ &3.20   &$6.17\times 10^{-6}$    &3.05&3.22	    &7.42	\\
80  &2  &$320\times 320$       &$1.09\times10^{-6}$ &3.01   &$7.69\times 10^{-7}$    &3.00&25.29	&7.84 	\\
160  &2  &$640\times 640$      &$1.36\times10^{-7}$ &3.00   &$9.60\times 10^{-8}$    &3.00&199.79	&7.90	\\
\hline
\end{tabular}
\caption{\footnotesize{Example 2, Linear scheme,
comparison of numerical errors and CPU times for computations on single-grid and sparse-grid. Lagrange interpolation for prolongation is used in sparse-grid computations.
Final time $T=1.0$. CFL number is $0.5$.
$N_r$: number of cells in each spatial direction of a root grid.
$N_L$: the finest level in a sparse-grid computation.
CPU: CPU time for a complete simulation. ``Ratio'' is the ratio of corresponding CPU times on an $N_h\times N_h$ mesh to that on a  $\frac{N_h}{2}\times \frac{N_h}{2}$ mesh. CPU time unit: seconds.}}
\label{tab:nonln_ln_2d}
\end{table}

\begin{table}[htbp]\footnotesize
\centering
\begin{tabular}{c c c c c c c c c}
\hline
\multicolumn{9}{c}{ Single-grid}\\ \hline
& & $N_h\times N_h$		&$L^\infty$ error		&Order  &$L^2$ error    &Order  &CPU(s)		&Ratio	\\
\hline
&&$80\times 80$         &$2.48\times10^{-4}$ &       &$1.05\times 10^{-4}$    &    &0.78	    &	\\
&&$160\times 160$       &$1.45\times10^{-5}$ &4.09   &$7.40\times 10^{-6}$    &3.83&6.17	    &7.90	\\
&&$320\times 320$       &$1.27\times10^{-6}$ &3.52   &$7.97\times 10^{-7}$    &3.21&48.85	    &7.92	\\
&&$640\times 640$       &$1.42\times10^{-7}$ &3.17   &$9.68\times 10^{-8}$    &3.04&399.74	&8.18	\\
\hline
\multicolumn{9}{c}{ Sparse-grid, refine root grids, $N_L=3$}\\ \hline
$N_r$   &$N_L$  &$N_h\times N_h$	&$L^\infty$ error   &Order	    &$L^2$ error    &Order  &CPU(s)	     &Ratio\\
\hline
10  &3  &$80\times 80$         &$2.28\times10^{-2}$ &       &$1.07\times 10^{-2}$    &    &0.58	    &	\\
20  &3  &$160\times 160$       &$9.45\times10^{-4}$ &4.59   &$4.36\times 10^{-4}$    &4.61&4.43	    &7.58	\\
40  &3  &$320\times 320$       &$1.17\times10^{-5}$ &6.33   &$4.44\times 10^{-6}$    &6.62&34.10	&7.70 	\\
80  &3  &$640\times 640$       &$1.85\times10^{-7}$ &5.98   &$9.93\times 10^{-8}$    &5.48&268.65	&7.88	\\
\hline
\multicolumn{9}{c}{ Sparse-grid, refine root grids, $N_L=2$}\\ \hline
$N_r$   &$N_L$  &$N_h\times N_h$	&$L^\infty$ error   &Order	    &$L^2$ error    &Order  &CPU(s)	     &Ratio\\
\hline
20  &2  &$80\times 80$         &$5.11\times10^{-3}$ &       &$2.32\times 10^{-3}$    &    &0.82	    &	\\
40  &2  &$160\times 160$       &$9.03\times10^{-5}$ &5.83   &$3.64\times 10^{-5}$    &5.99&6.29	    &7.65	\\
80  &2  &$320\times 320$       &$1.60\times10^{-6}$ &5.82   &$8.24\times 10^{-7}$    &5.47&49.78	&7.92 	\\
160 &2  &$640\times 640$       &$1.43\times10^{-7}$ &3.49   &$9.68\times 10^{-8}$    &3.09&403.93	&8.11	\\
\hline
\end{tabular}
\caption{\footnotesize{Example 2, WENO scheme,
comparison of numerical errors and CPU times for computations on single-grid and sparse-grid. Lagrange interpolation for prolongation is used in sparse-grid computations.
Final time $T=1.0$. CFL number is $0.5$.
$N_r$: number of cells in each spatial direction of a root grid.
$N_L$: the finest level in a sparse-grid computation.
CPU: CPU time for a complete simulation. ``Ratio'' is the ratio of corresponding CPU times on an $N_h\times N_h$ mesh to that on a  $\frac{N_h}{2}\times \frac{N_h}{2}$ mesh. CPU time unit: seconds.}}
\label{tab:nonln_weno_2d}
\end{table}

\begin{table}[htbp]\small
\centering
\begin{tabular}{c c c c c c c}
\hline
\multicolumn{7}{c}{ Single-grid}\\ \hline
&&$N_h\times N_h$		&CPU(s)	&	&Grid point \# &	\\
\hline
&&$80\times 80$         &0.78	 &   &6,561	&\\
&&$160\times 160$       &6.17	 &   &25,921 &	\\
&&$320\times 320$       &48.85	&&103,041	&\\
&&$640\times 640$       &399.74	&&410,881	&\\
\hline
\multicolumn{7}{c}{ Sparse-grid, refine root grids, $N_L=3$}\\ \hline
$N_r$   &$N_L$  &$N_h\times N_h$	&CPU(s)	     &sparse/single (C) &Grid point \#  &sparse/single (G)\\
\hline
10  &3  &$80\times 80$         &0.58	    &0.75  &4,847	&0.74\\
20  &3  &$160\times 160$       &4.43	    &0.72  &18,487	&0.71\\
40  &3  &$320\times 320$       &34.10	    &0.70  &72,167 	&0.70\\
80  &3  &$640\times 640$       &268.65	    &0.67  &285,127	&0.69\\
\hline
\multicolumn{7}{c}{ Sparse-grid, refine root grids, $N_L=2$}\\ \hline
$N_r$   &$N_L$  &$N_h\times N_h$	&CPU(s)	     &sparse/single (C) &Grid point \#  &sparse/single (G)\\
\hline
20  &2  &$80\times 80$         &0.82	    &1.05  &6,805	&1.04\\
40  &2  &$160\times 160$       &6.29	    &1.02  &26,405	&1.02\\
80  &2  &$320\times 320$       &49.78	    &1.02  &104,005 &1.01\\
160 &2  &$640\times 640$       &403.93	    &1.01  &412,805	&1.00\\
\hline
\end{tabular}
\caption{\footnotesize{Example 2, comparison of numbers of spatial grid points and CPU times in sparse-grid and single-grid computations. $N_r$: number of cells in each spatial direction of a root grid. $N_L$: the finest level in a sparse-grid computation.
``CPU": CPU times in WENO computations (Table \ref{tab:nonln_weno_2d}, unit: seconds); ``sparse/single (C)": ratios of corresponding CPU times on sparse grids to that on single grids; ``Grid point \#": numbers of spatial grid points in sparse-grid and single-grid computations; ``sparse/single (G)": ratios of corresponding numbers of spatial grid points on sparse grids to that on single grids.}}
\label{tab:gridcount_2d}
\end{table}

\paragraph{Example 3 (A 3D Nonlinear equation):}
\begin{equation}
\label{eqn:Ex3}
\begin{cases}
u_t+(\frac{1}{2}u^2)_x+(\frac{1}{2}u^2)_y+(\frac{1}{2}u^2)_z=-0.1u, \qquad 0 \leq x \leq 4\pi, 0 \leq y \leq 4\pi, 0 \leq z \leq 4\pi;\\
u(x,y,z,0)=\sin(x-0.5y-0.5z),
\end{cases}
\end{equation}
with periodic boundary condition. This is a 3D version of Example 2. The exact solution of the problem is $u(x,y,z,t)=e^{-0.1t}\sin(x-0.5y-0.5z)$.
We compute this 3D nonlinear problem till final time $T=1$ by both the third order linear scheme and WENO scheme
on both single grids and sparse grids, and compare their computational efficiency.
In Table \ref{tab:nonln_ln_3d} and Table \ref{tab:nonln_weno_3d},
 the $L^\infty$ errors, $L^2$ errors, the corresponding numerical accuracy orders, and CPU times on
successively refined meshes to show the efficiency of computations on sparse grids are presented.
Similar as the previous examples, for the linear scheme, comparable numerical errors and third order accuracy are obtained on single grids and sparse grids,
except that the errors are larger on coarser meshes for the $N_L=3$ case.
For the WENO scheme, numerical errors and accuracy on sparse grids are comparable to those on single grids if the mesh is more refined, while
on relatively coarse meshes, the numerical errors of computations on sparse grids are larger than that on single grids.
Again, the computations on sparse grids show superconvergence to reach comparable numerical errors and accuracy as that on single grids.
For CPU times, we observe about $70\% \sim 85\%$ saving if $N_L=3$, and $50\% \sim 65\%$ saving for the $N_L=2$ case.
We also count the saving of grid points used in the sparse-grid computations.
For 3D sparse grids with $N_r$ number of cells in each spatial direction of the root grid and the finest level $N_L$, the number of
grid points of the most refined grid $\Omega^{N_L,N_L,N_L}$ is
 $(N_r \cdot 2^{N_L}+1)^3$. This is also the number of
grid points which are used in discretizing and solving the PDEs on the corresponding 3D single grids. While in
the 3D sparse-grid computations, the number of
grid points which are used in discretizing and solving the PDEs is
\begin{eqnarray}
\label{sparsegridnum3D}
&\displaystyle\sum_{i=0}^{N_L}\sum_{j=0}^{N_L-i}(N_r \cdot 2^i+1) \cdot (N_r \cdot 2^j+1) \cdot (N_r \cdot 2^{N_L-i-j}+1) \nonumber\\
&\displaystyle+\sum_{i=0}^{N_L-1}\sum_{j=0}^{N_L-i-1}(N_r \cdot 2^i+1) \cdot (N_r \cdot 2^j+1) \cdot (N_r \cdot 2^{N_L-i-j-1}+1)  \nonumber\\
&\displaystyle+\sum_{i=0}^{N_L-2}\sum_{j=0}^{N_L-i-2}(N_r \cdot 2^i+1) \cdot (N_r \cdot 2^j+1) \cdot (N_r \cdot 2^{N_L-i-j-2}+1).
\end{eqnarray}
Again the count of grid points on sparse grids does not include those grid points used in the prolongation step, since the prolongation step is done only at the final time step and does not directly involve solving the PDEs.
In Table \ref{tab:gridcount_3d}, we list numbers of spatial grid points which are used in discretizing and solving this 3D example on sparse grids and single grids.
It is observed that sparse-grid computations just use about $20\%$ grid points of single grids for the $N_L=3$ case, and about $50\%$ grid points of single grids for the $N_L=2$ case.
The ratios of numbers of spatial grid points on sparse grids to that on single grids are compared with the ratios of CPU times
 on sparse grids to that on single grids. Consistent results are obtained and verify that the saving of CPU times in sparse-grid computations is due to the fact that less number of grid points is used.

\begin{table}[htbp]\footnotesize
\centering
\begin{tabular}{c c c c c c c c }
\hline
\multicolumn{8}{c}{ Single-grid}\\ \hline
& & $N_h\times N_h\times N_h$	&$L^\infty$ error		&Order	&$L^2$ error &Order &CPU(s)	\\
\hline
&&$80\times 80\times 80$         &$3.82\times10^{-4}$ &       &$2.38\times 10^{-4}$    &    &31.81	    	\\
&&$160\times 160\times 160$       &$4.20\times10^{-5}$ &3.19   &$2.82\times 10^{-5}$    &3.08&447.45    	\\
&&$320\times 320\times 320$       &$4.99\times10^{-6}$ &3.07   &$3.47\times 10^{-6}$    &3.02&10,618.80 	\\
&&$640\times 640\times 640$       &$6.14\times10^{-7}$ &3.02  &$4.33\times 10^{-7}$   &3.01&160,662.44  	\\
\hline
\multicolumn{8}{c}{ Sparse-grid, refine root grids, $N_L=3$}\\ \hline
$N_r$   &$N_L$  &$N_h\times N_h\times N_h$	&$L^\infty$ error   &Order &$L^2$ error   &Order	&CPU(s)	     \\
\hline
10  &3  &$80\times 80\times 80$         &$7.09\times10^{-2}$ &       &$1.38\times 10^{-2}$    &    &13.89	    	\\
20  &3  &$160\times 160\times 160$       &$1.10\times10^{-4}$ &9.33   &$4.50\times 10^{-5}$    &8.26&115.09       	\\
40  &3  &$320\times 320\times 320$       &$5.12\times10^{-6}$ &4.43   &$3.52\times 10^{-6}$    &3.68&1,641.56	    	\\
80  &3  &$640\times 640\times 640$       &$6.16\times10^{-7}$ &3.06   &$4.33\times 10^{-7}$    &3.02&26,387.20	\\
\hline
\multicolumn{8}{c}{ Sparse-grid, refine root grids, $N_L=2$}\\ \hline
$N_r$   &$N_L$  &$N_h\times N_h\times N_h$	&$L^\infty$ error   &Order &$L^2$ error    &Order	&CPU(s)	     \\
\hline
20  &2  &$80\times 80\times 80$         &$5.03\times10^{-4}$ &       &$2.85\times 10^{-4}$    &    &14.83	    	\\
40  &2  &$160\times 160\times 160$       &$4.29\times10^{-5}$ &3.55   &$2.85\times 10^{-5}$    &3.32&218.43       	\\
80  &2  &$320\times 320\times 320$       &$5.00\times10^{-6}$ &3.10   &$3.48\times 10^{-6}$    &3.03&3,563.63	     	\\
160  &2  &$640\times 640\times 640$      &$6.14\times10^{-7}$ &3.02   &$4.33\times 10^{-7}$    &3.01&75,802.70    	\\
\hline
\end{tabular}
\caption{\footnotesize{Example 3, Linear scheme,
comparison of numerical errors and CPU times for computations on single-grid and sparse-grid. Lagrange interpolation for prolongation is used in sparse-grid computations.
Final time $T=1.0$. CFL number is $0.75$.
$N_r$: number of cells in each spatial direction of a root grid.
$N_L$: the finest level in a sparse-grid computation.
CPU: CPU time for a complete simulation. CPU time unit: seconds.}}
\label{tab:nonln_ln_3d}
\end{table}

\begin{table}[htbp]\small
\centering
\begin{tabular}{c c c c c c c}
\hline
\multicolumn{7}{c}{ Single-grid}\\ \hline
&&$N_h\times N_h\times N_h$		&CPU(s)	&	&Grid point \# &	\\
\hline
&&$80\times 80\times 80$          &31.81	    &&531,441	&\\
&&$160\times 160\times 160$       &447.45	    &&4,173,281	 &\\
&&$320\times 320\times 320$       &10,618.80	&&33,076,161	&\\
&&$640\times 640\times 640$       &160,662.44    &&263,374,721  &\\
\hline
\multicolumn{7}{c}{ Sparse-grid, refine root grids, $N_L=3$}\\ \hline
$N_r$   &$N_L$  &$N_h\times N_h\times N_h$	&CPU(s)	     &sparse/single (C) &Grid point \#  &sparse/single (G)\\
\hline
10  &3  &$80\times 80\times 80$         &13.89	    &0.44         &132,549    &0.25\\
20  &3  &$160\times 160\times 160$       &115.09	    &0.26     &967,679    &0.23\\
40  &3  &$320\times 320\times 320$       &1,641.56	    &0.15     &7,385,739  &0.22\\
80  &3  &$640\times 640\times 640$       &26,387.20	    &0.16        &57,693,059  &0.22\\
\hline
\multicolumn{7}{c}{ Sparse-grid, refine root grids, $N_L=2$}\\ \hline
$N_r$   &$N_L$  &$N_h\times N_h\times N_h$	&CPU(s)	     &sparse/single (C) &Grid point \#  &sparse/single (G)\\
\hline
20  &2  &$80\times 80\times 80$         &14.83	    &0.47         &276,570    &0.52\\
40  &2  &$160\times 160\times 160$      &218.43	    &0.49         &2,096,330    &0.50\\
80  &2  &$320\times 320\times 320$      &3,563.63	    &0.34     &16,317,450  &0.49\\
160  &2  &$640\times 640\times 640$     &75,802.70	    &0.47       &128,750,090  &0.49\\
\hline
\end{tabular}
\caption{\footnotesize{Example 3, comparison of numbers of spatial grid points and CPU times in sparse-grid and single-grid computations. $N_r$: number of cells in each spatial direction of a root grid. $N_L$: the finest level in a sparse-grid computation.
``CPU": CPU times in computations by Linear scheme (Table \ref{tab:nonln_ln_3d}, unit: seconds); ``sparse/single (C)": ratios of corresponding CPU times on sparse grids to that on single grids; ``Grid point \#": numbers of spatial grid points in sparse-grid and single-grid computations; ``sparse/single (G)": ratios of corresponding numbers of spatial grid points on sparse grids to that on single grids.}}
\label{tab:gridcount_3d}
\end{table}

\begin{table}[htbp]\footnotesize
\centering
\begin{tabular}{c c c c c c c c }
\hline
\multicolumn{8}{c}{ Single-grid}\\ \hline
& & $N_h\times N_h\times N_h$	&$L^\infty$ error		&Order	&$L^2$ error &Order &CPU(s)	\\
\hline
&&$80\times 80\times 80$         &$2.73\times10^{-3}$ &       &$1.11\times 10^{-3}$    &    &54.20	    	\\
&&$160\times 160\times 160$       &$1.34\times10^{-4}$ &4.34   &$5.71\times 10^{-5}$    &4.29&816.97    	\\
&&$320\times 320\times 320$       &$7.95\times10^{-6}$ &4.08   &$4.10\times 10^{-6}$    &3.80&16,298.10 	\\
&&$640\times 640\times 640$       &$7.08\times10^{-7}$ &3.49  &$4.47\times 10^{-7}$   &3.20   &217,742.62   	\\
\hline
\multicolumn{8}{c}{ Sparse-grid, refine root grids, $N_L=3$}\\ \hline
$N_r$   &$N_L$  &$N_h\times N_h\times N_h$	&$L^\infty$ error   &Order &$L^2$ error   &Order	&CPU(s)	     \\
\hline
10  &3  &$80\times 80\times 80$         &$1.86\times10^{-1}$ &       &$4.05\times 10^{-2}$    &    &14.68	    	\\
20  &3  &$160\times 160\times 160$       &$2.19\times10^{-2}$ &3.09   &$5.62\times 10^{-3}$    &2.85&197.14       	\\
40  &3  &$320\times 320\times 320$       &$5.69\times10^{-4}$ &5.26   &$1.44\times 10^{-4}$    &5.28&3,005.57	     	\\
80  &3  &$640\times 640\times 640$       &$1.64\times10^{-6}$ &8.44   &$6.15\times 10^{-7}$    &7.87&48,099.20		\\
\hline
\multicolumn{8}{c}{ Sparse-grid, refine root grids, $N_L=2$}\\ \hline
$N_r$   &$N_L$  &$N_h\times N_h\times N_h$	&$L^\infty$ error   &Order &$L^2$ error    &Order	&CPU(s)	     \\
\hline
20  &2  &$80\times 80\times 80$         &$3.27\times10^{-2}$ &       &$1.03\times 10^{-2}$    &    &27.05	    	\\
40  &2  &$160\times 160\times 160$       &$1.29\times10^{-3}$ &4.67   &$3.75\times 10^{-4}$    &4.77&408.71     	\\
80  &2  &$320\times 320\times 320$       &$1.45\times10^{-5}$ &6.47   &$5.52\times 10^{-6}$    &6.09&6,385.60	 	\\
160  &2  &$640\times 640\times 640$      &$7.32\times10^{-7}$ &4.31   &$4.49\times 10^{-7}$    &3.62&115,599.00 	\\
\hline
\end{tabular}
\caption{\footnotesize{Example 3, WENO scheme,
comparison of numerical errors and CPU times for computations on single-grid and sparse-grid. Lagrange interpolation for prolongation is used in sparse-grid computations.
Final time $T=1.0$. CFL number is $0.75$.
$N_r$: number of cells in each spatial direction of a root grid.
$N_L$: the finest level in a sparse-grid computation.
CPU: CPU time for a complete simulation. CPU time unit: seconds.}}
\label{tab:nonln_weno_3d}
\end{table}

\subsection{Examples with Discontinuous Solutions}

In this section, we test the method for solving 2D and 3D Burgers' equations, in which shock waves form at certain time and discontinuities develop in the solutions.
WENO prolongation in sparse grid combination techniques is necessary to deal with discontinuous solutions here.

\begin{figure}[p]
    \centering
    \begin{subfigure}[b]{0.48\textwidth}
        \includegraphics[width=0.8\textwidth]{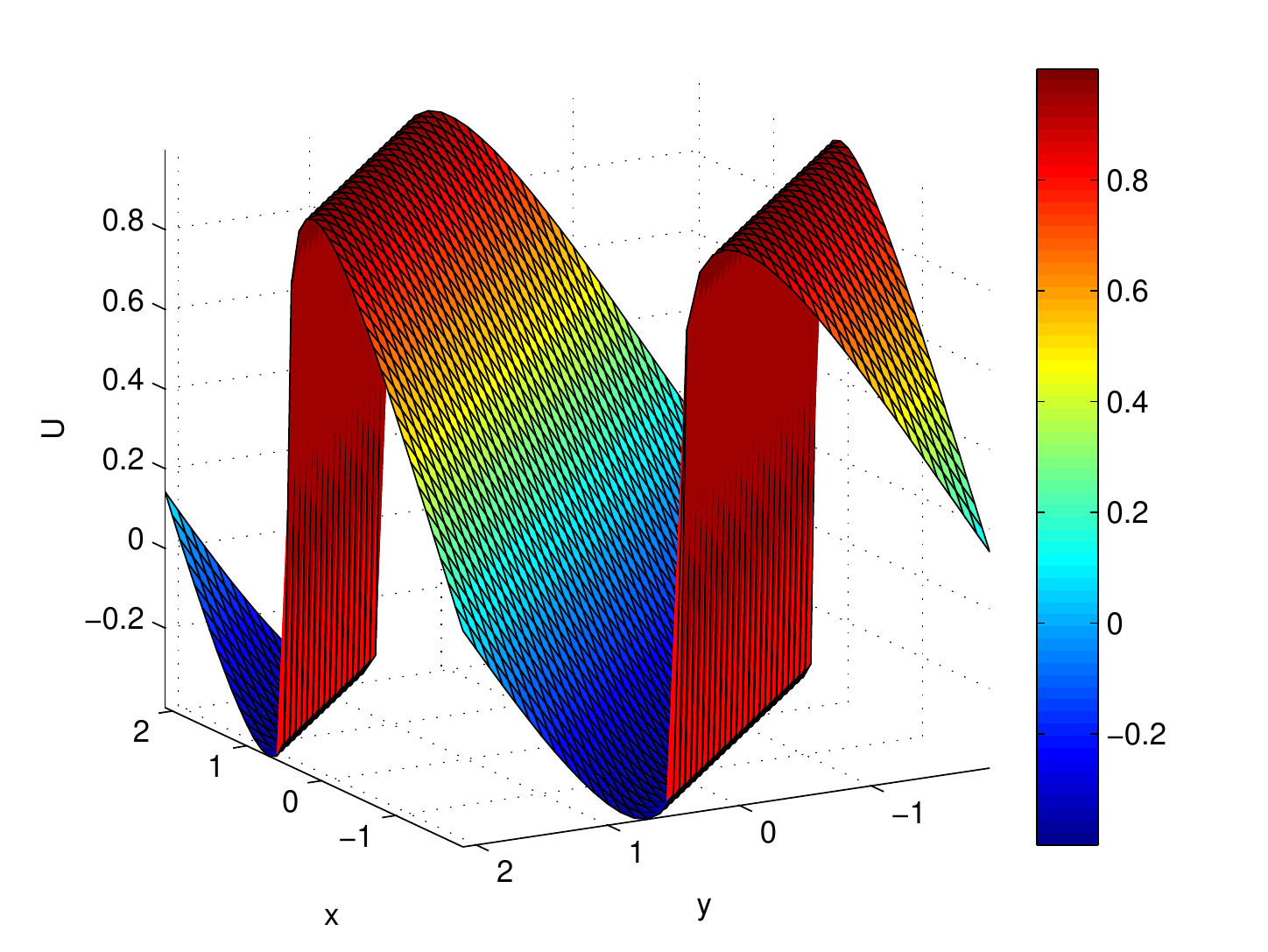}
        \caption{}
        \label{40_3_surf}
    \end{subfigure}
        \begin{subfigure}[b]{0.48\textwidth}
        \includegraphics[width=0.8\textwidth]{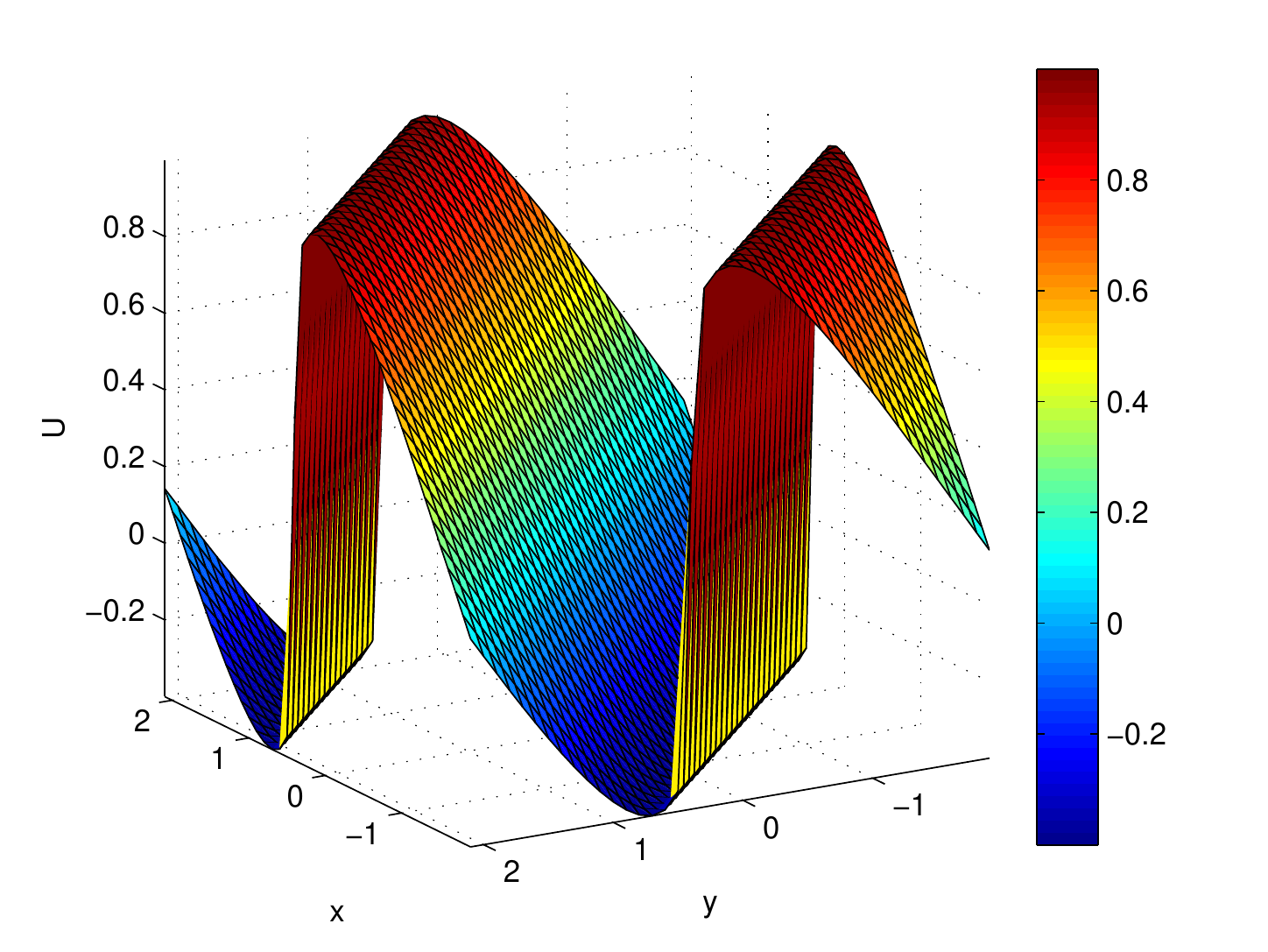}
        \caption{}
        \label{40_3_surf_single}
    \end{subfigure}
        \begin{subfigure}[b]{0.48\textwidth}
        \includegraphics[width=0.8\textwidth]{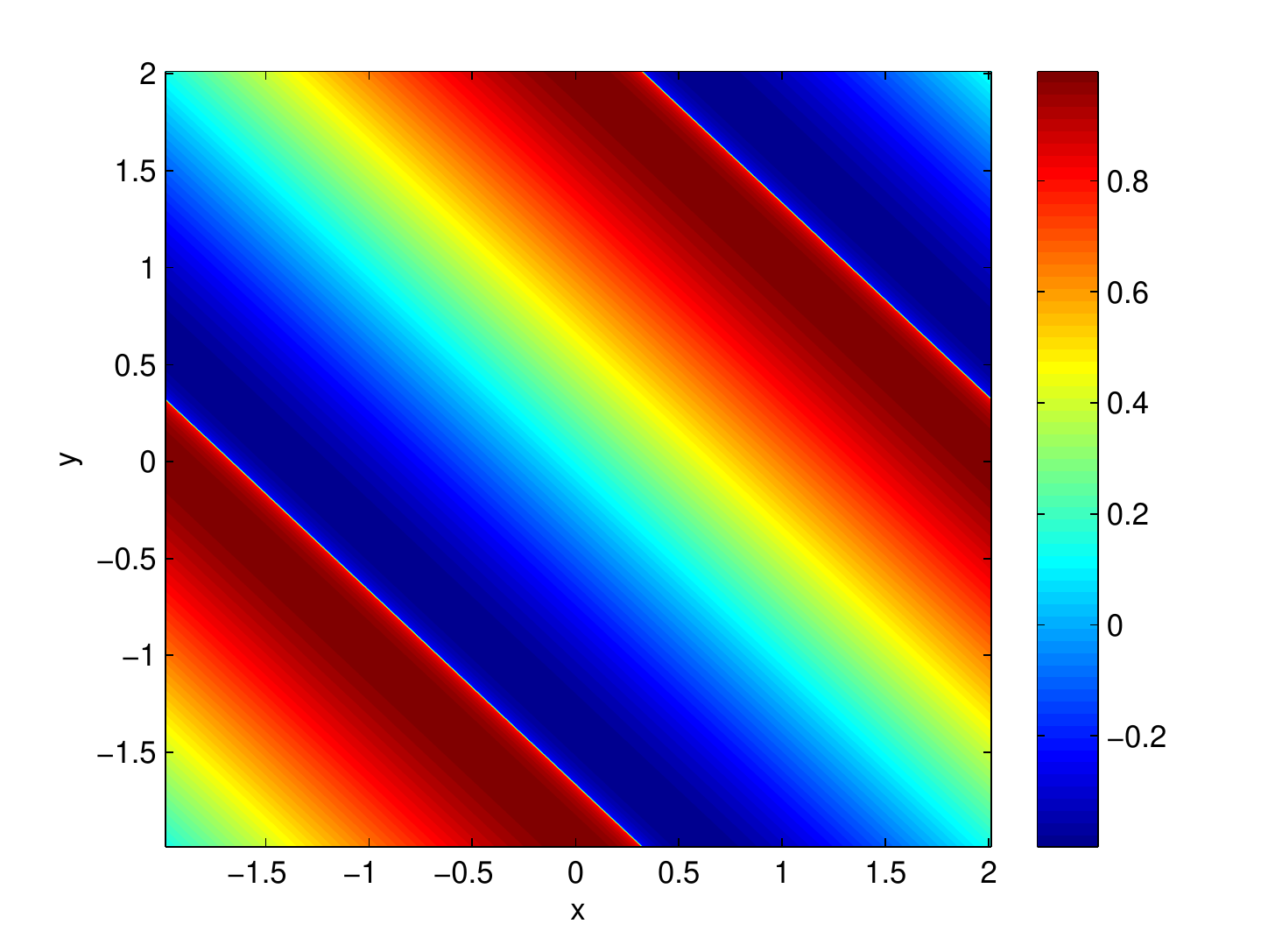}
        \caption{}
        \label{40_3_contour}
    \end{subfigure}
        \begin{subfigure}[b]{0.48\textwidth}
        \includegraphics[width=0.8\textwidth]{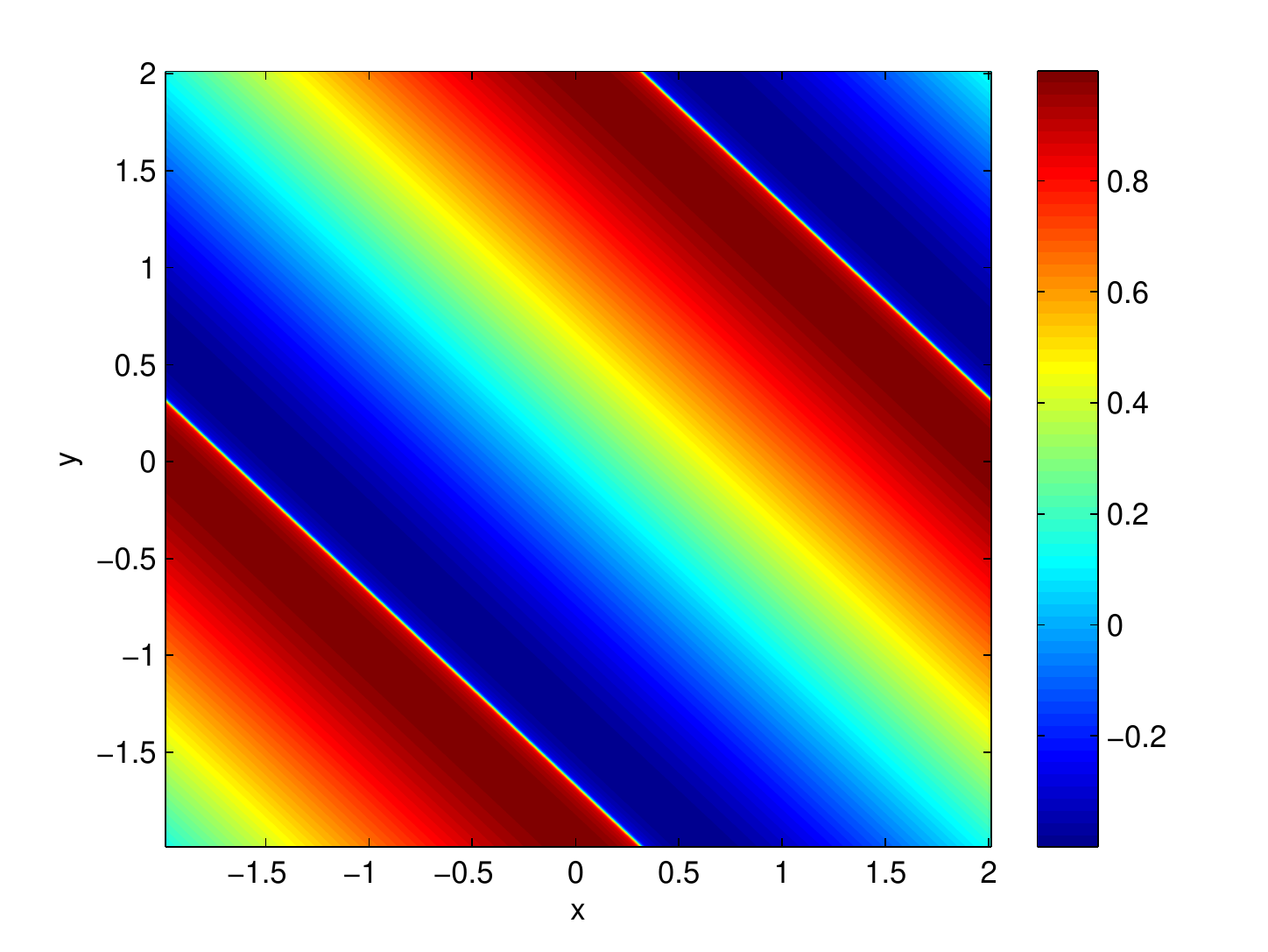}
        \caption{}
        \label{40_3_contour_single}
    \end{subfigure}
        \begin{subfigure}[b]{0.48\textwidth}
        \includegraphics[width=0.8\textwidth]{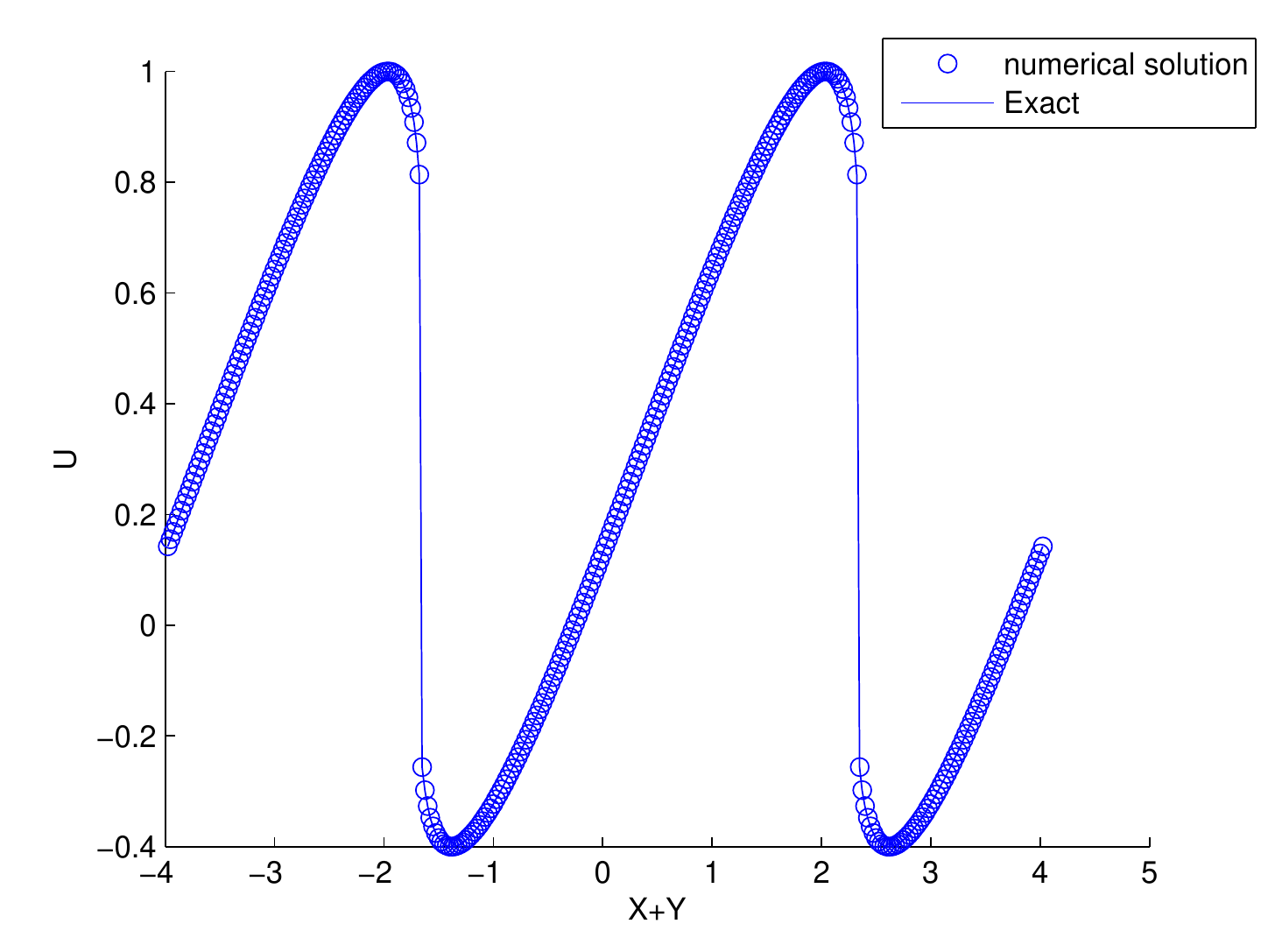}
        \caption{}
        \label{40_3_curve}
    \end{subfigure}
        \begin{subfigure}[b]{0.48\textwidth}
        \includegraphics[width=0.8\textwidth]{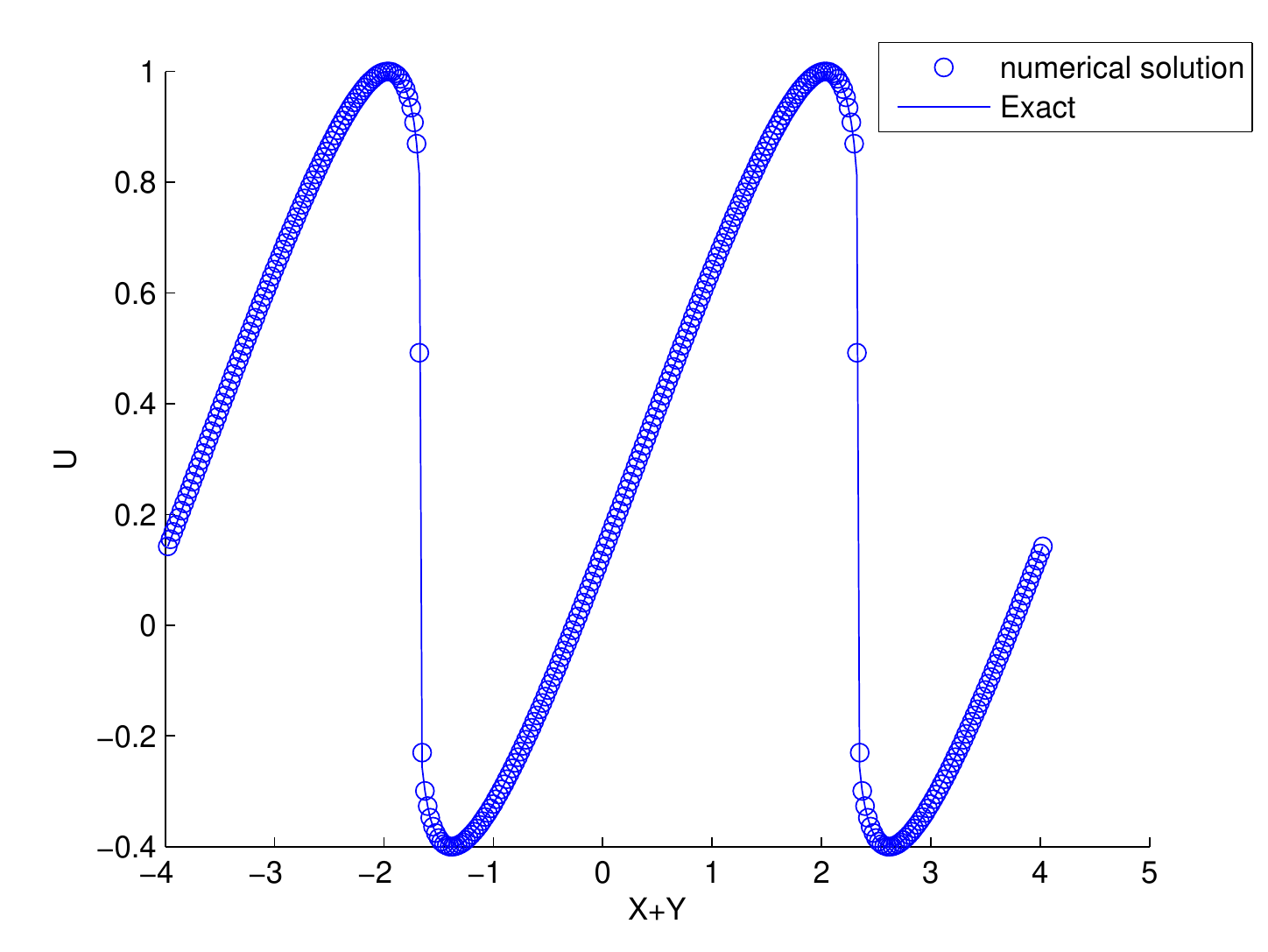}
        \caption{}
        \label{40_3_curve_single}
    \end{subfigure}
    \caption{Example 4, solution of two dimensional Burgers' equation by third order WENO scheme on sparse grids ($N_r=40$ for root grid, finest level $N_L=3$ in the sparse-grid computation) and the corresponding $320 \times 320$ single grid, using third order WENO interpolation for prolongation in sparse grid combination. $CFL=0.5$, final time $T=5/\pi^2$. (a), (c), (e): sparse-grid results; (b), (d), (f): single-grid results. (a), (b): 3D surface plots of the solutions; (c), (d): contour plots; (e), (f): 1D cutting-plot along $x=y$. }
\label{2dplot_burg}
\end{figure}

\paragraph{Example 4 (A 2D Burgers' equation):}
\begin{equation}
\label{eqn:Ex4}
\begin{cases}
u_t+(\frac{u^2}{2})_x+(\frac{u^2}{2})_y=0, \qquad (x,y)\in [-2,2] \times [-2,2];\\
u(x,y,0)=0.3+0.7\sin(\frac{\pi}{2}(x+y)),
\end{cases}
\end{equation}
with periodic boundary conditions. We first apply both the third order linear scheme and WENO scheme
on single grids and sparse grids to solve this problem to $T=0.5/\pi^2$, when the solution is still smooth.
The WENO prolongation in sparse grid combination techniques is used in the WENO scheme. In Table \ref{tab:bg_ln_2d}
and Table \ref{tab:bg_weno_2d}, the $L^\infty$ errors, $L^2$ errors, the corresponding numerical accuracy orders, and CPU times on
successively refined meshes are reported. We observe that on relatively coarse meshes, the sparse-grid computations have larger errors than
the single-grid computations. However along with the mesh refinement, the numerical errors of the sparse-grid computations catch up with that
of the single-grid computations, and comparable numerical errors are obtained. As that in the previous examples, the WENO scheme show superconvergence when the mesh is refined. Then we compute the solution of the problem at $T=5/\pi^2$, when the shock waves form and the solution is discontinuous. The results on sparse grids and the corresponding single grid
are shown in Figure \ref{2dplot_burg}. We can see that the numerical solution by the sparse grid WENO scheme with the WENO prolongation is similar as that by the single-grid computation.
The non-oscillatory and high resolution properties of the WENO scheme for resolving shock waves are preserved well in the sparse-grid computations. In Table \ref{ex9_cpu_time}, CPU time costs for the sparse-grid and single-grid computations are listed. We observe about $25\% \sim 35\%$ CPU time savings by using
 sparse-grid with $N_L=3$ for this 2D problem with discontinuous solution.

\begin{table}[htbp]\footnotesize
\centering
\begin{tabular}{c c c c c c c c}
\hline
\multicolumn{8}{c}{ Single-grid}\\ \hline
& & $N_h\times N_h$		&$L^\infty$ error		&Order  &$L^2$ error    &Order  &CPU(s)			\\
\hline
&&$80\times 80$         &$6.76\times10^{-6}$ &       &$3.88\times 10^{-6}$    &    &0.04	    	\\
&&$160\times 160$       &$7.84\times10^{-7}$ &3.11   &$4.67\times 10^{-7}$    &3.05&0.25	    	\\
&&$320\times 320$       &$9.78\times10^{-8}$ &3.00   &$5.73\times 10^{-8}$    &3.03&1.92	    	\\
&&$640\times 640$       &$1.22\times10^{-8}$ &3.00   &$7.10\times 10^{-9}$    &3.01&15.13 		\\
\hline
\multicolumn{8}{c}{ Sparse-grid, refine root grids, $N_L=3$}\\ \hline
$N_r$   &$N_L$  &$N_h\times N_h$	&$L^\infty$ error   &Order	    &$L^2$ error    &Order  &CPU(s)	     \\
\hline
10  &3  &$80\times 80$         &$9.29\times10^{-5}$ &       &$3.01\times 10^{-5}$    &    &0.03	    	\\
20  &3  &$160\times 160$       &$1.04\times10^{-5}$ &3.16   &$3.94\times 10^{-6}$    &2.93&0.16	    	\\
40  &3  &$320\times 320$       &$1.34\times10^{-6}$ &2.96   &$5.03\times 10^{-7}$    &2.97&1.00	     	\\
80  &3  &$640\times 640$       &$1.67\times10^{-7}$ &3.00   &$6.30\times 10^{-8}$    &3.00&7.22	    	\\
\hline
\multicolumn{8}{c}{ Sparse-grid, refine root grids, $N_L=2$}\\ \hline
$N_r$   &$N_L$  &$N_h\times N_h$	&$L^\infty$ error   &Order	    &$L^2$ error    &Order  &CPU(s)	     \\
\hline
20  &2  &$80\times 80$         &$1.57\times10^{-5}$ &       &$5.08\times 10^{-6}$    &    &0.03	    	\\
40  &2  &$160\times 160$       &$1.59\times10^{-6}$ &3.31   &$6.05\times 10^{-7}$    &3.07&0.21	    	\\
80  &2  &$320\times 320$       &$1.91\times10^{-7}$ &3.05   &$7.44\times 10^{-8}$    &3.02&1.51	     	\\
160 &2  &$640\times 640$       &$2.36\times10^{-8}$ &3.02   &$9.20\times 10^{-9}$    &3.01&11.50		\\
\hline
\end{tabular}
\caption{\footnotesize{Example 4, Linear scheme,
comparison of numerical errors and CPU times for computations on single-grid and sparse-grid. Lagrange interpolation for prolongation is used in sparse-grid computations.
Final time $T=0.5/\pi^2$. CFL number is $0.5$.
$N_r$: number of cells in each spatial direction of a root grid.
$N_L$: the finest level in a sparse-grid computation.
CPU: CPU time for a complete simulation. CPU time unit: seconds.}}
\label{tab:bg_ln_2d}
\end{table}

\begin{table}[htbp]\footnotesize
\centering
\begin{tabular}{c c c c c c c c}
\hline
\multicolumn{8}{c}{ Single-grid}\\ \hline
& & $N_h\times N_h$		&$L^\infty$ error		&Order  &$L^2$ error    &Order  &CPU(s)			\\
\hline
&&$80\times 80$         &$7.80\times 10^{-5}$ &       &$2.09\times 10^{-6}$    &    &0.08	    	\\
&&$160\times 160$       &$3.05\times 10^{-6}$ &4.68   &$9.90\times 10^{-7}$    &4.40&0.57	    	\\
&&$320\times 320$       &$1.66\times 10^{-7}$ &4.20   &$6.99\times 10^{-8}$    &3.82&4.29	    	\\
&&$640\times 640$       &$1.42\times 10^{-8}$ &3.55   &$7.42\times 10^{-9}$    &3.24&34.11 		\\
\hline
\multicolumn{8}{c}{ Sparse-grid, refine root grids, $N_L=3$}\\ \hline
$N_r$   &$N_L$  &$N_h\times N_h$	&$L^\infty$ error   &Order	    &$L^2$ error    &Order  &CPU(s)	     \\
\hline
10  &3  &$80\times 80$         &$3.40\times10^{-2}$ &       &$3.82\times 10^{-3}$    &    &0.08	    	\\
20  &3  &$160\times 160$       &$1.72\times10^{-3}$ &4.30   &$1.72\times 10^{-4}$    &4.47&0.49	    	\\
40  &3  &$320\times 320$       &$1.02\times10^{-5}$ &7.39   &$1.03\times 10^{-6}$    &7.39&3.30	     	\\
80  &3  &$640\times 640$       &$2.78\times10^{-8}$ &8.52   &$7.79\times 10^{-9}$    &7.04&23.92	    	\\
\hline
\multicolumn{8}{c}{ Sparse-grid, refine root grids, $N_L=2$}\\ \hline
$N_r$   &$N_L$  &$N_h\times N_h$	&$L^\infty$ error   &Order	    &$L^2$ error    &Order  &CPU(s)	     \\
\hline
20  &2  &$80\times 80$         &$1.02\times10^{-2}$ &       &$1.13\times 10^{-3}$    &    &0.10	    	\\
40  &2  &$160\times 160$       &$6.76\times10^{-5}$ &7.23   &$9.96\times 10^{-6}$    &6.82&0.63	    	\\
80  &2  &$320\times 320$       &$2.77\times10^{-7}$ &7.93   &$7.33\times 10^{-8}$    &7.09&4.49	     	\\
160 &2  &$640\times 640$       &$1.44\times10^{-8}$ &4.27   &$7.43\times 10^{-9}$    &3.30&33.90		\\
\hline
\end{tabular}
\caption{\footnotesize{Example 4, WENO scheme,
comparison of numerical errors and CPU times for computations on single-grid and sparse-grid. WENO interpolation for prolongation is used in sparse-grid computations.
Final time $T=0.5/\pi^2$. CFL number is $0.5$.
$N_r$: number of cells in each spatial direction of a root grid.
$N_L$: the finest level in a sparse-grid computation.
CPU: CPU time for a complete simulation. CPU time unit: seconds.}}
\label{tab:bg_weno_2d}
\end{table}

\begin{table}[htbp]
\centering
\begin{tabular}{c c c}
\hline
$N_h\times N_h$  &CPU time on sparse-grid    &CPU time on single-grid\\
\hline
$80\times 80$       &0.50   &0.64\\
$160\times 160$     &3.24    &5.11\\
$320\times 320$     &30.08    &40.55\\
$640\times 640$     &225.01    &341.36\\
\hline
\end{tabular}
\caption{\footnotesize{Example 4, WENO scheme, comparison of CPU times for computations of discontinuous solution on single-grid and sparse-grid. WENO interpolation for prolongation is used in
sparse-grid computations. Final time $T=5/\pi^2$. CFL number is $0.5$. $N_L=3$ in sparse-grid computations. CPU time unit: seconds.}}
\label{ex9_cpu_time}
\end{table}

\paragraph{Example 5 (A 3D Burgers' equation):}
\begin{equation}
\label{eqn:Ex5}
\begin{cases}
u_t+(\frac{u^2}{2})_x+(\frac{u^2}{2})_y+(\frac{u^2}{2})_z=0, \qquad (x,y,z)\in [-3,3] \times [-3,3] \times [-3,3]; \\
u(x,y,z,0)=0.3+0.7\sin(\frac{\pi}{2}(x+y+z)),
\end{cases}
\end{equation}
with periodic boundary conditions. As that for the last example, we first apply both the third order linear scheme and WENO scheme
on single grids and sparse grids to solve this problem to $T=0.5/\pi^2$, when the solution is still smooth. Both the regular Lagrange prolongation and the
WENO prolongation in sparse grid combination techniques are tested in the WENO scheme. In Table \ref{tab:bg_ln_3d}
 and Table \ref{tab:bg_weno_3d}, the $L^\infty$ errors, $L^2$ errors, the corresponding numerical accuracy orders, and CPU times on
successively refined meshes are reported. For the linear scheme, comparable numerical errors and third order accuracy are obtained on single grids and sparse grids,
except that the errors are larger on the $80 \times 80 \times 80$ mesh for the $N_L=3$ case. For the WENO scheme, numerical errors and accuracy on sparse grids are comparable to those on single grids if the mesh is more refined, while
on relatively coarse meshes, the numerical errors of computations on sparse grids are larger than that on single grids, especially for the $N_L=3$ case. The $N_L=2$ case has
smaller numerical errors and better accuracy, but more CPU time costs than these of the $N_L=3$ case. The WENO computations on sparse grids show superconvergence to reach comparable numerical errors and accuracy as that on single grids. The numerical errors of the sparse-grid computations using the WENO prolongation are larger than that using
the regular Lagrange prolongation on relatively coarse meshes, however they are comparable if the mesh is more refined. This is due to the typical asymptotic convergence property of nonlinear WENO approximations.
For the 3D problem, much more significant CPU time savings than the 2D problem are observed. We observe up to $77\%$ CPU time saving by using
 sparse-grid with $N_L=3$, and up to $56\%$ CPU time saving if $N_L=2$ in the simulations.

\begin{figure}[p]
    \centering
    \begin{subfigure}[b]{0.48\textwidth}
        \includegraphics[width=0.8\textwidth]{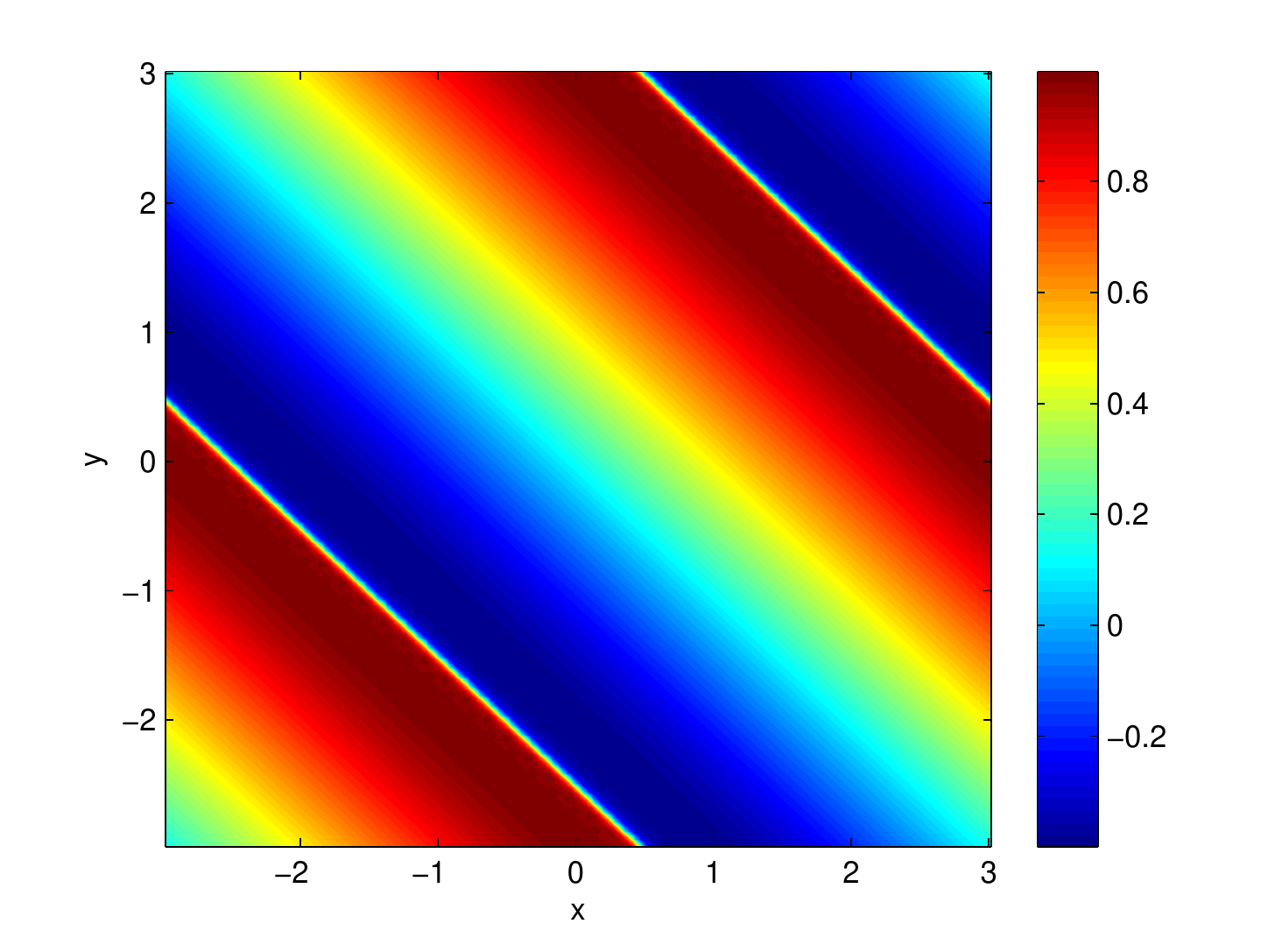}
        \caption{}
        \label{3d_40_3_contour_0}
    \end{subfigure}
       \begin{subfigure}[b]{0.48\textwidth}
        \includegraphics[width=0.8\textwidth]{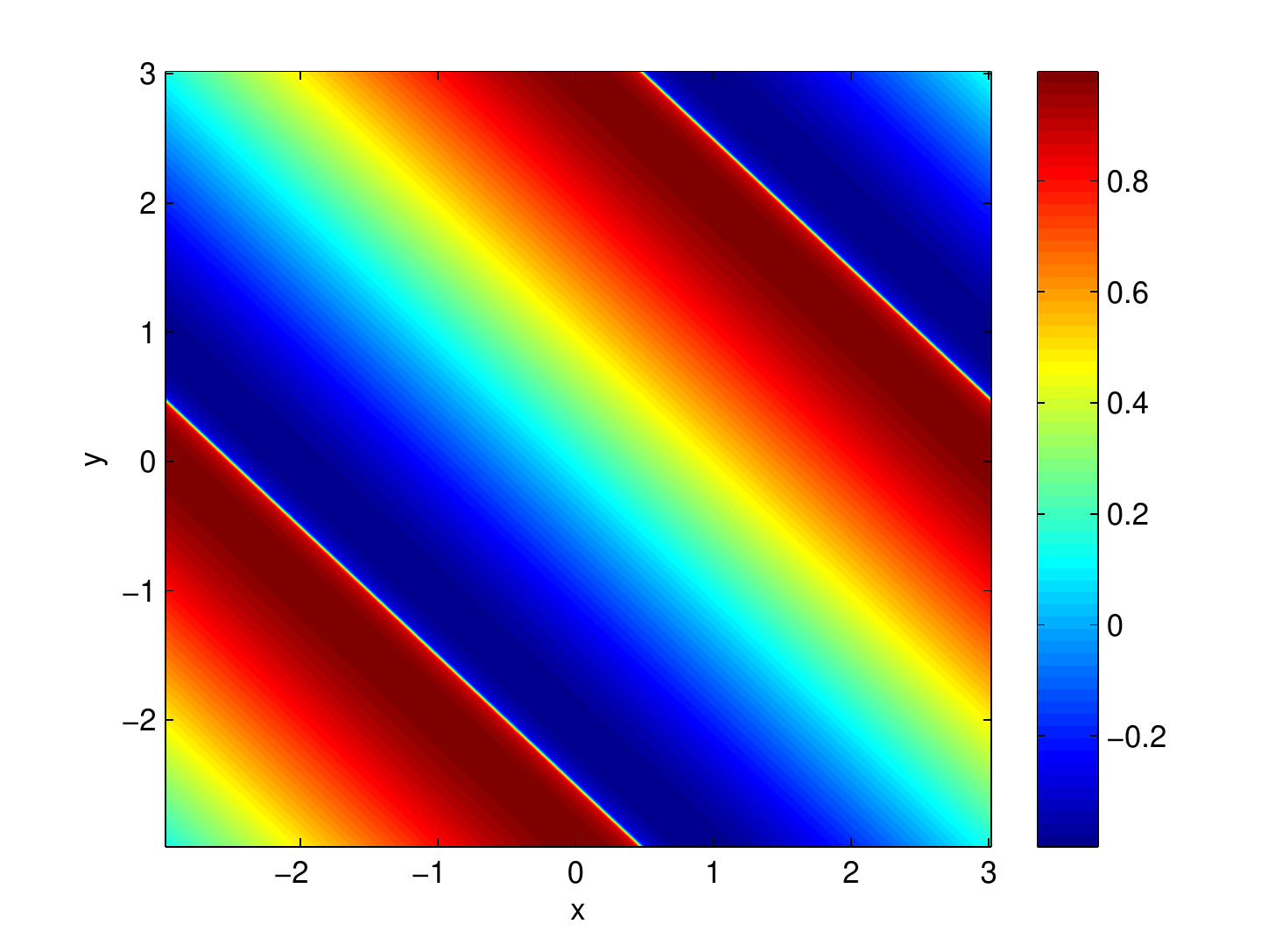}
        \caption{}
        \label{3d_40_3_contour_0_single}
    \end{subfigure}
        \begin{subfigure}[b]{0.48\textwidth}
        \includegraphics[width=0.8\textwidth]{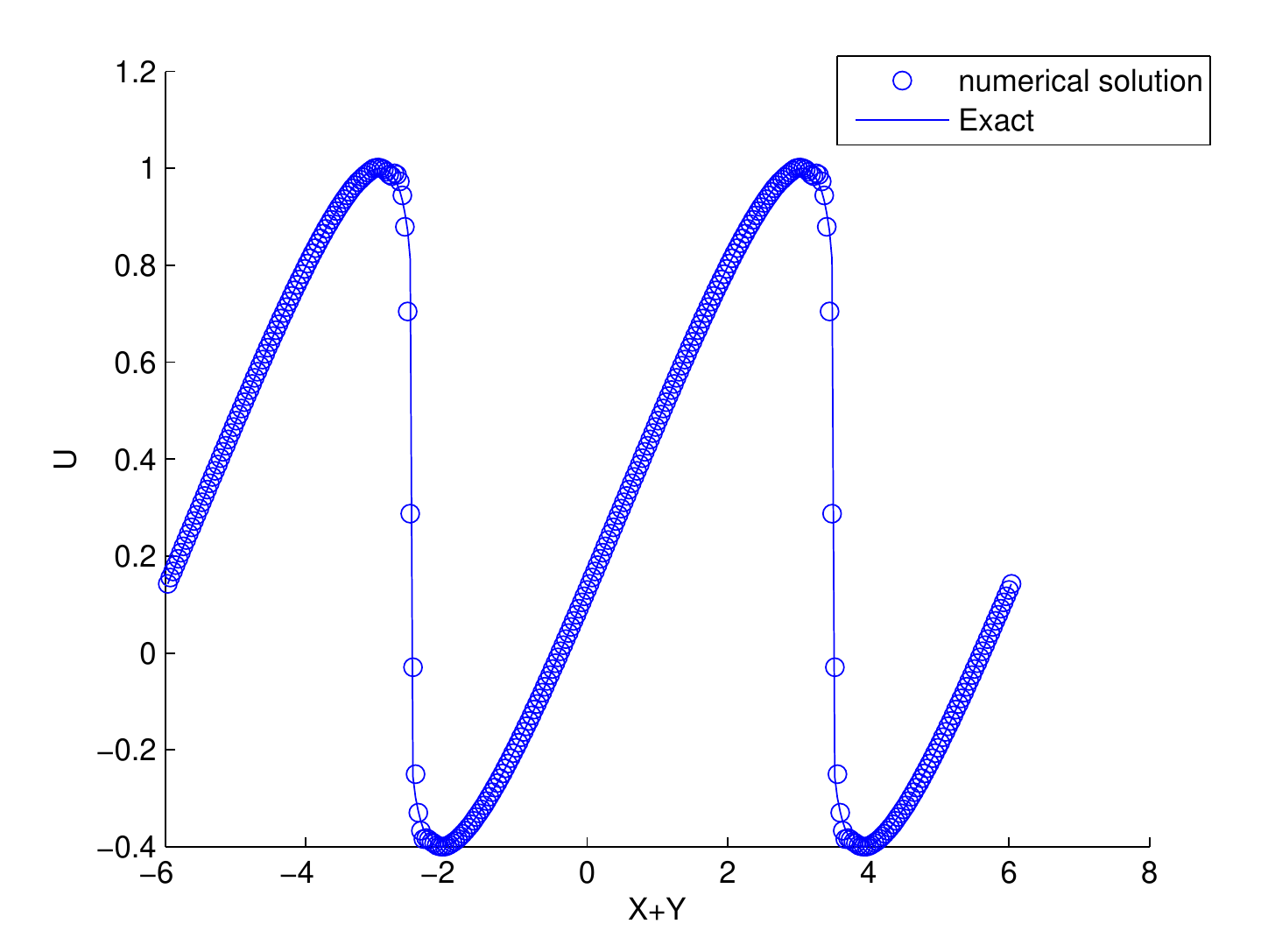}
        \caption{}
        \label{3d_40_3_curve_0}
    \end{subfigure}
        \begin{subfigure}[b]{0.48\textwidth}
        \includegraphics[width=0.8\textwidth]{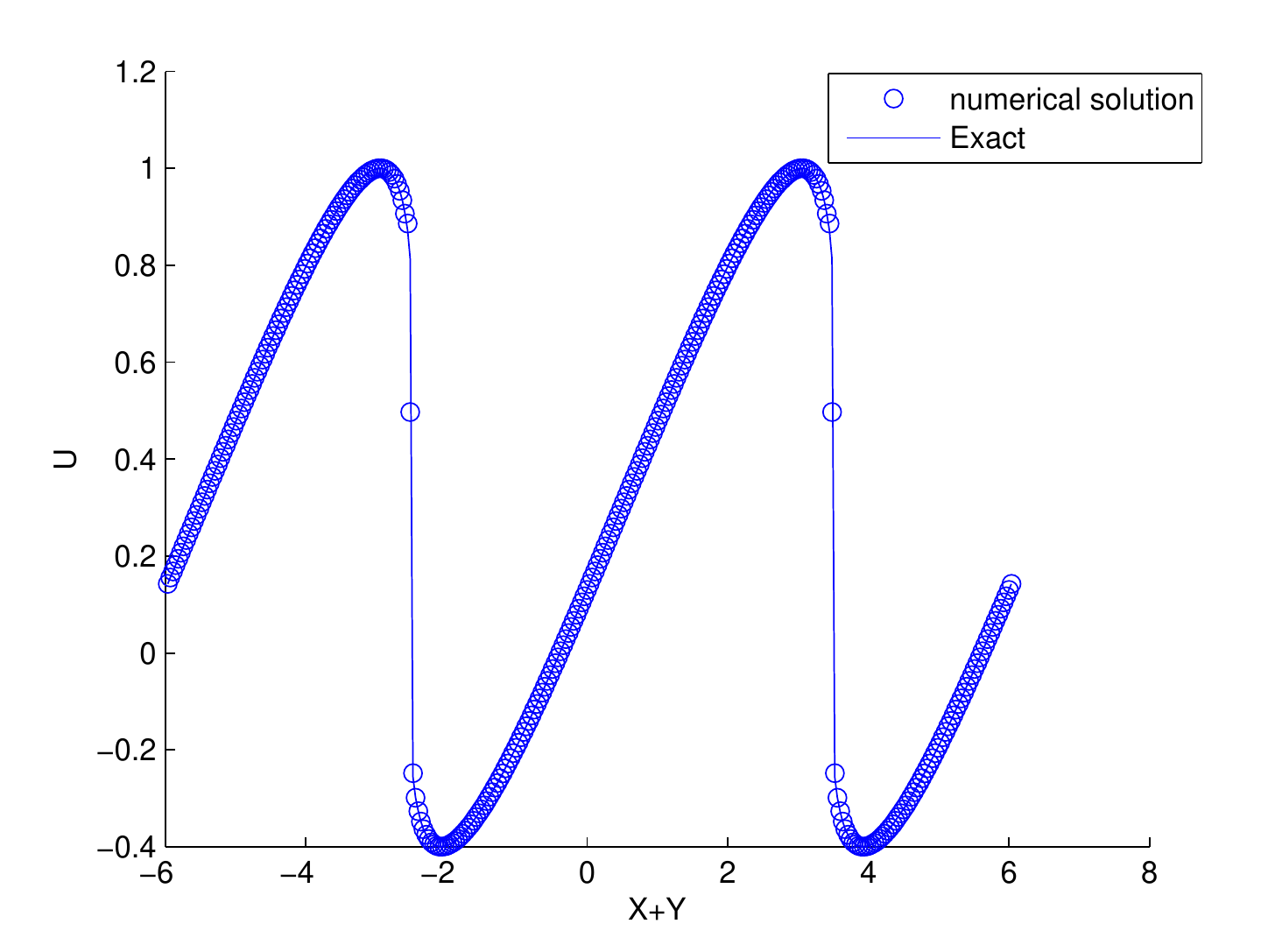}
        \caption{}
        \label{3d_40_3_curve_0_single}
    \end{subfigure}
        \begin{subfigure}[b]{0.48\textwidth}
        \includegraphics[width=0.8\textwidth]{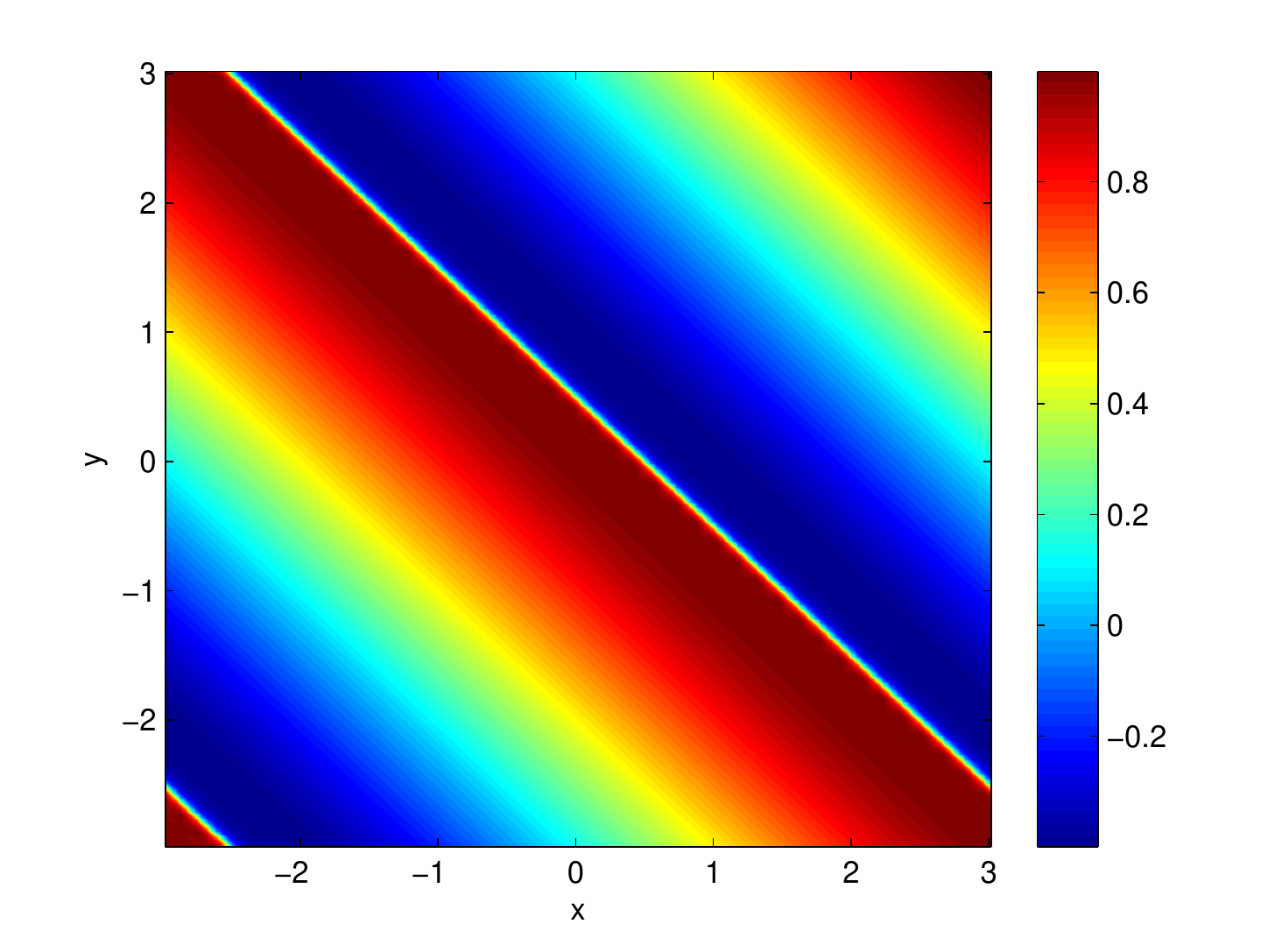}
        \caption{}
        \label{3d_40_3_contour_160}
    \end{subfigure}
        \begin{subfigure}[b]{0.48\textwidth}
        \includegraphics[width=0.8\textwidth]{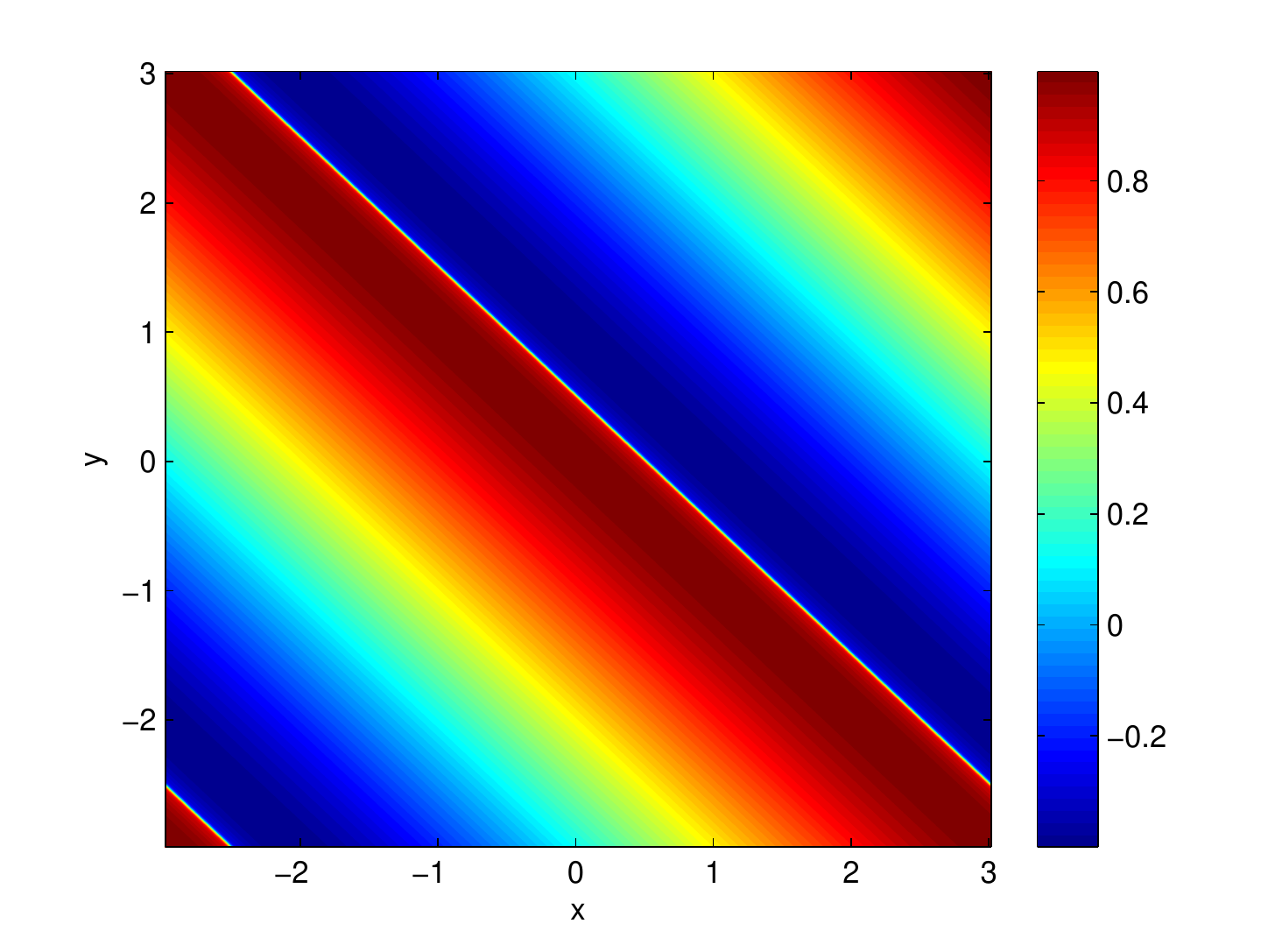}
        \caption{}
        \label{3d_40_3_contour_160_single}
    \end{subfigure}
        \begin{subfigure}[b]{0.48\textwidth}
        \includegraphics[width=0.8\textwidth]{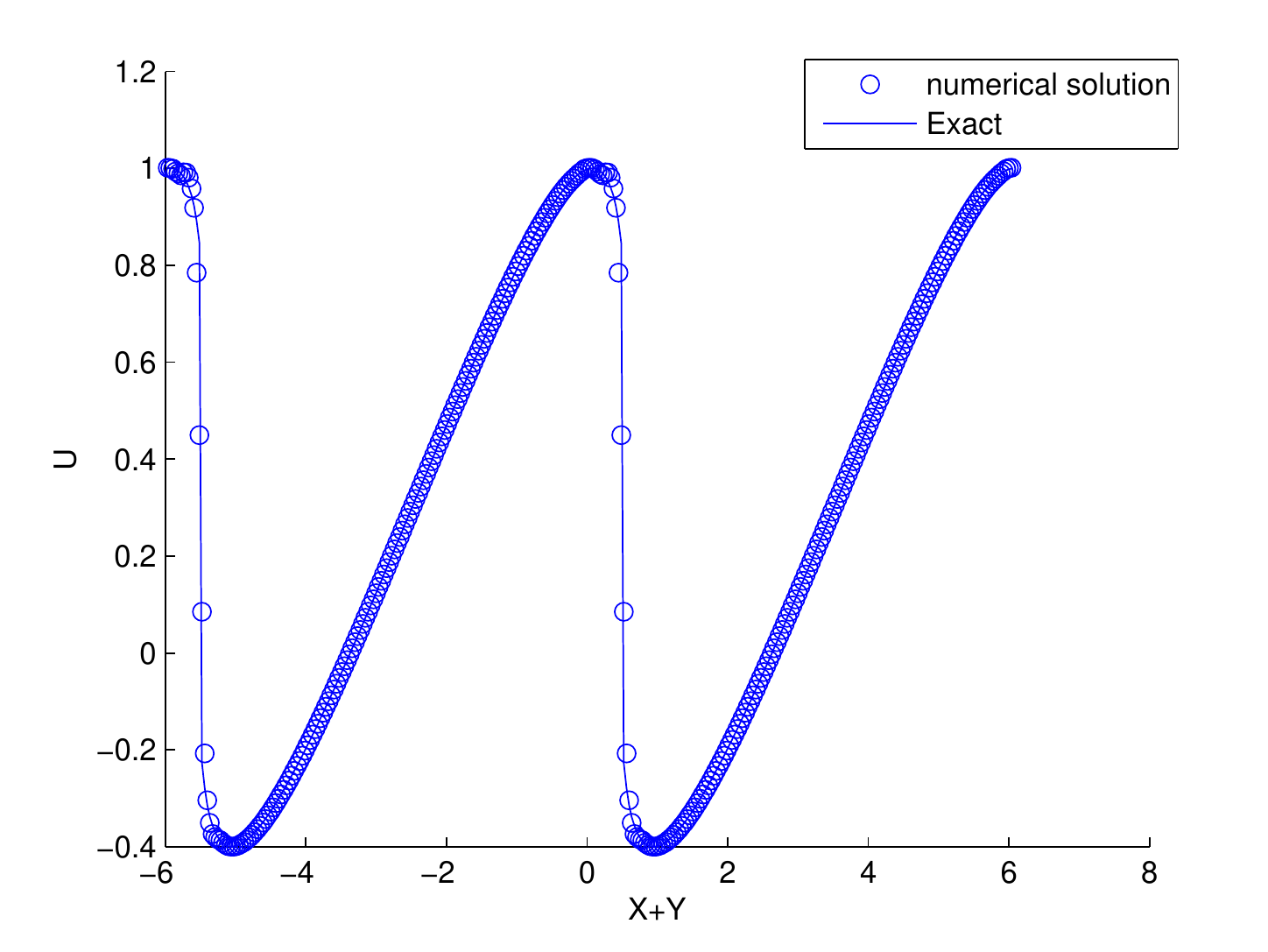}
        \caption{}
        \label{3d_40_3_curve_160}
    \end{subfigure}
        \begin{subfigure}[b]{0.48\textwidth}
        \includegraphics[width=0.8\textwidth]{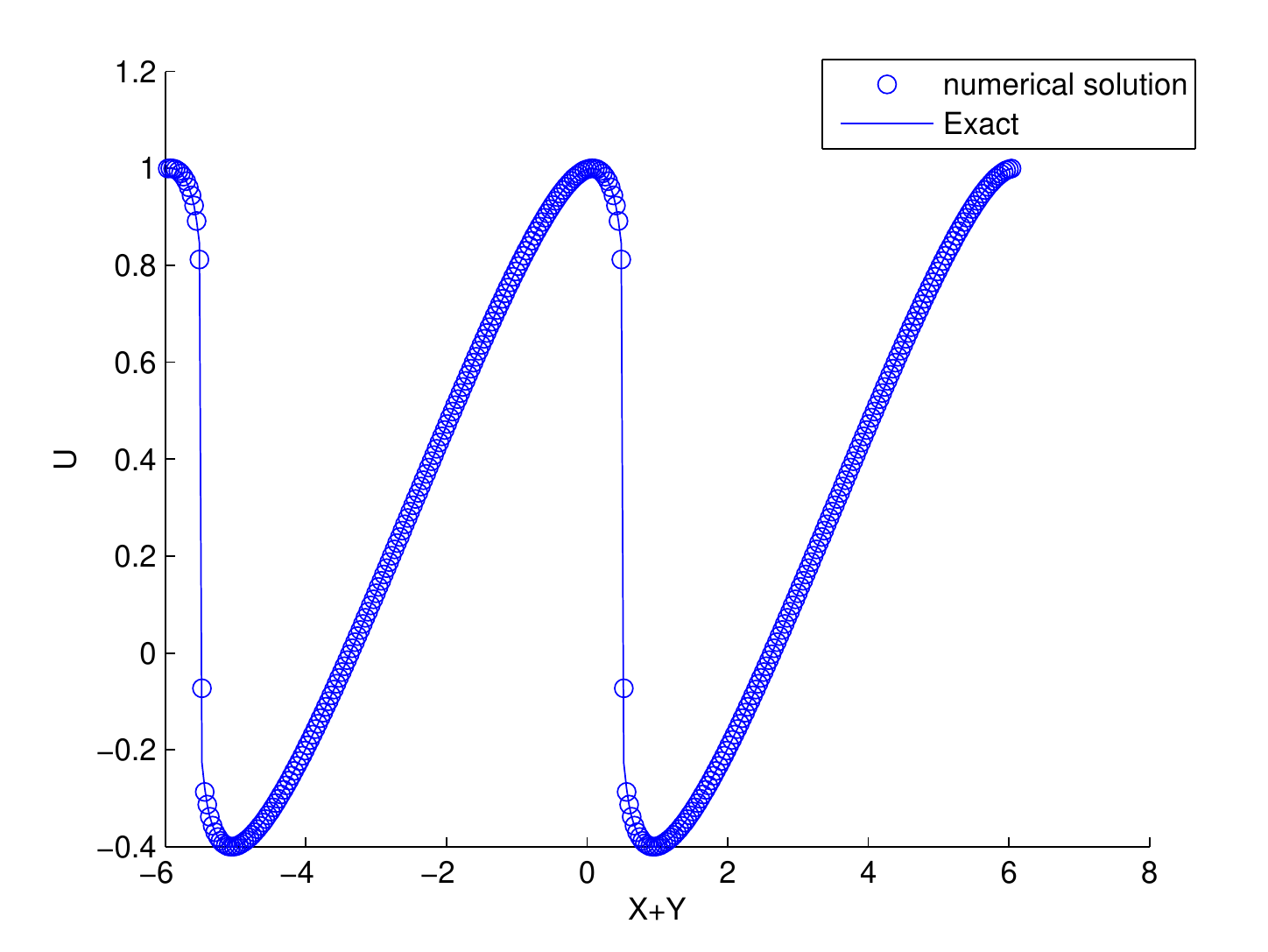}
        \caption{}
        \label{3d_40_3_curve_160_single}
    \end{subfigure}
    \caption{Example 5, solution of three dimensional Burgers' equation by third order WENO scheme on sparse grids ($N_r=40$ for root grid, finest level $N_L=3$ in the sparse-grid computation) and the corresponding $320 \times 320\times 320$ single grid, using third order WENO interpolation for prolongation in sparse grid combination. $CFL=0.75$, final time $T=5/\pi^2$.  (a), (c), (e),(g): sparse-grid results; (b), (d), (f),(h): single-grid results. (a), (b): 2D contour plots of $x-y$ plane cutting at $z=-3$; (c), (d): 1D cutting-plot along $x=y$ on the plane $z=-3$; (e), (f): 2D contour plots of $x-y$ plane cutting at $z=0$; (g), (h): 1D cutting-plot along $x=y$ on the plane $z=0$.
    }
    \label{3dplot_burg}
\end{figure}

Then we compute the solution of the problem at $T=5/\pi^2$, when the shock waves form and the solution is discontinuous. If we apply the WENO scheme on sparse grids using the  regular Lagrange prolongation, spurious oscillations are observed in the results. These spurious oscillations are removed if the WENO prolongation is used in sparse grid combinations. We show the results on sparse grids and the corresponding single grid
in Figure \ref{3dplot_burg}. Similar as the 2D example, we can see that the numerical solution by the sparse grid WENO scheme with the WENO prolongation is similar as that by the single-grid computation.
The non-oscillatory and high resolution properties of the WENO scheme for resolving shock waves are preserved well in the sparse-grid computations. In Table \ref{ex10_cpu_time}, CPU time costs for the sparse-grid and single-grid computations are listed. We observe about $50\% \sim 70\%$ CPU time savings by using
 sparse-grid with $N_L=3$ for this 3D problem with discontinuous solution.

\begin{table}[htbp]\footnotesize
\centering
\begin{tabular}{c c c c c c c c }
\hline
\multicolumn{8}{c}{ Single-grid}\\ \hline
& & $N_h\times N_h\times N_h$	&$L^\infty$ error		&Order	&$L^2$ error &Order &CPU(s)	\\
\hline
&&$80\times 80\times 80$          &$8.31\times10^{-6}$ &       &$4.08\times 10^{-6}$    &    &3.59	   \\
&&$160\times 160\times 160$       &$1.13\times10^{-6}$ &2.88   &$5.05\times 10^{-7}$    &3.01&57.95       	\\
&&$320\times 320\times 320$       &$1.45\times10^{-7}$ &2.96   &$6.28\times 10^{-8}$    &3.01&1,072.84    	\\
&&$640\times 640\times 640$       &$1.83\times10^{-8}$ &2.98   &$7.83\times 10^{-9}$    &3.00&16,671.93  	\\
\hline
\multicolumn{8}{c}{ Sparse-grid, refine root grids, $N_L=3$}\\ \hline
$N_r$   &$N_L$  &$N_h\times N_h\times N_h$	&$L^\infty$ error   &Order &$L^2$ error   &Order	&CPU(s)	     \\
\hline
10  &3  &$80\times 80\times 80$          &$9.74\times10^{-5}$ &       &$1.44\times 10^{-5}$    &    &3.64	    	\\
20  &3  &$160\times 160\times 160$       &$1.29\times10^{-6}$ &6.24   &$5.91\times 10^{-7}$    &4.61&37.98      	\\
40  &3  &$320\times 320\times 320$       &$1.49\times10^{-7}$ &3.12   &$6.54\times 10^{-8}$    &3.17&344.57	     	\\
80  &3  &$640\times 640\times 640$       &$1.84\times10^{-8}$ &3.01   &$7.91\times 10^{-9}$    &3.05&4,314.76		\\
\hline
\multicolumn{8}{c}{ Sparse-grid, refine root grids, $N_L=2$}\\ \hline
$N_r$   &$N_L$  &$N_h\times N_h\times N_h$	&$L^\infty$ error   &Order &$L^2$ error    &Order	&CPU(s)	     \\
\hline
20  &2  &$80\times 80\times 80$         &$8.69\times10^{-6}$ &       &$3.82\times 10^{-6}$    &     &2.96	    	\\
40  &2  &$160\times 160\times 160$       &$1.12\times10^{-6}$ &2.96   &$4.98\times 10^{-7}$    &2.94&35.67      	\\
80  &2  &$320\times 320\times 320$       &$1.45\times10^{-7}$ &2.95   &$6.26\times 10^{-8}$    &2.99&473.95	     	\\
160  &2  &$640\times 640\times 640$      &$1.83\times10^{-8}$ &2.98   &$7.83\times 10^{-9}$    &3.00&8,552.87   	\\
\hline
\end{tabular}
\caption{\footnotesize{Example 5, Linear scheme,
comparison of numerical errors and CPU times for computations on single-grid and sparse-grid. Lagrange interpolation for prolongation is used in sparse-grid computations.
Final time $T=0.5/\pi^2$. CFL number is $0.75$.
$N_r$: number of cells in each spatial direction of a root grid.
$N_L$: the finest level in a sparse-grid computation.
CPU: CPU time for a complete simulation. CPU time unit: seconds.}}
\label{tab:bg_ln_3d}
\end{table}

\begin{table}[htbp]\footnotesize
\centering
\begin{tabular}{c c c c c c c c }
\hline
\multicolumn{8}{c}{ Single-grid}\\ \hline
& & $N_h\times N_h\times N_h$	&$L^\infty$ error		&Order	&$L^2$ error &Order &CPU(s)	\\
\hline
&&$80\times 80\times 80$         &$1.57\times 10^{-5}$ &       &$5.71\times 10^{-6}$    &    &6.38	    	\\
&&$160\times 160\times 160$       &$1.36\times 10^{-6}$ &3.53   &$5.54\times 10^{-7}$    &3.37&109.56	    	\\
&&$320\times 320\times 320$       &$1.52\times 10^{-7}$ &3.16   &$6.43\times 10^{-8}$    &3.11&1,839.92		\\
&&$640\times 640\times 640$       &$1.85\times 10^{-8}$ &3.04   &$7.88\times 10^{-9}$    &3.03&28,522.89		\\
\hline
\multicolumn{8}{c}{ Sparse-grid, refine root grids, $N_L=3$, Lagrange interpolation for prolongation}\\ \hline
$N_r$   &$N_L$  &$N_h\times N_h\times N_h$	&$L^\infty$ error   &Order &$L^2$ error   &Order	&CPU(s)	     \\
\hline
10  &3  &$80\times 80\times 80$          &$1.39\times10^{-2}$ &       &$1.37\times 10^{-3}$    &    &4.28	    	\\
20  &3  &$160\times 160\times 160$       &$3.17\times10^{-4}$ &5.45   &$2.12\times 10^{-5}$    &6.01&42.97      	\\
40  &3  &$320\times 320\times 320$       &$2.82\times10^{-7}$ &10.14  &$7.51\times 10^{-8}$    &8.14&488.31	     	\\
80  &3  &$640\times 640\times 640$       &$1.90\times10^{-8}$ &3.89   &$8.01\times 10^{-9}$    &3.23&6,566.05		\\
\hline
\multicolumn{8}{c}{ Sparse-grid, refine root grids, $N_L=2$, Lagrange interpolation for prolongation}\\ \hline
$N_r$   &$N_L$  &$N_h\times N_h\times N_h$	&$L^\infty$ error   &Order &$L^2$ error    &Order	&CPU(s)	     \\
\hline
20  &2  &$80\times 80\times 80$         &$7.49\times10^{-4}$ &        &$9.36\times 10^{-5}$    &     &5.06	    	\\
40  &2  &$160\times 160\times 160$       &$2.06\times10^{-6}$ &8.51   &$5.56\times 10^{-7}$    &7.40&59.90      	\\
80  &2  &$320\times 320\times 320$       &$1.52\times10^{-7}$ &3.76   &$6.41\times 10^{-8}$    &3.12&802.99	     	\\
160  &2  &$640\times 640\times 640$      &$1.85\times10^{-8}$ &3.04   &$7.87\times 10^{-9}$    &3.02&14,057.20  	\\
\hline
\multicolumn{8}{c}{ Sparse-grid, refine root grids, $N_L=3$, WENO interpolation for prolongation}\\ \hline
$N_r$   &$N_L$  &$N_h\times N_h\times N_h$	&$L^\infty$ error   &Order &$L^2$ error   &Order	&CPU(s)	     \\
\hline
10  &3  &$80\times 80\times 80$         &$2.14\times10^{-1}$ &       &$1.12\times 10^{-2}$    &    &6.68	    	\\
20  &3  &$160\times 160\times 160$      &$9.30\times10^{-4}$ &7.85   &$8.26\times 10^{-5}$    &7.08&62.41	    	\\
40  &3  &$320\times 320\times 320$      &$6.01\times10^{-7}$ &10.59  &$9.44\times 10^{-8}$   &9.77&641.35	     	\\
80  &3  &$640\times 640\times 640$      &$1.95\times10^{-8}$ &4.95   &$8.00\times 10^{-9}$    &3.56&7,696.12		\\
\hline
\multicolumn{8}{c}{ Sparse-grid, refine root grids, $N_L=2$, WENO interpolation for prolongation}\\ \hline
$N_r$   &$N_L$  &$N_h\times N_h\times N_h$	&$L^\infty$ error   &Order &$L^2$ error   &Order	&CPU(s)	     \\
\hline
20  &2  &$80\times 80\times 80$          &$4.48\times10^{-3}$ &       &$4.51\times 10^{-4}$    &    &6.05	    	\\
40  &2  &$160\times 160\times 160$       &$3.79\times10^{-6}$ &10.21  &$7.55\times 10^{-7}$    &9.22&73.17	    	\\
80  &2  &$320\times 320\times 320$       &$1.56\times10^{-7}$ &4.60   &$6.39\times 10^{-8}$    &3.56&916.54	     	\\
160 &2  &$640\times 640\times 640$       &$1.85\times10^{-8}$ &3.07   &$7.87\times 10^{-9}$    &3.02&15,179.40		\\
\hline
\end{tabular}
\caption{\footnotesize{Example 5, WENO scheme,
comparison of numerical errors and CPU times for computations on single-grid and sparse-grid. Both Lagrange and WENO interpolations for prolongation are used in sparse-grid computations.
Final time $T=0.5/\pi^2$. CFL number is $0.75$.
$N_r$: number of cells in each spatial direction of a root grid.
$N_L$: the finest level in a sparse-grid computation.
CPU: CPU time for a complete simulation. CPU time unit: seconds.}}
\label{tab:bg_weno_3d}
\end{table}

\begin{table}[htbp]
\centering
\begin{tabular}{c c c}
\hline
$N_h\times N_h\times N_h$  &CPU time on sparse-grid    &CPU time on single-grid\\
\hline
$80\times 80\times 80$        &19.18       &28.54\\
$160\times 160\times 160$     &241.01      &484.41\\
$320\times 320\times 320$     &3,518.58     &8,023.08\\
$640\times 640\times 640$     &53,515.70     &172,030.00\\
\hline
\end{tabular}
\caption{\footnotesize{Example 5, WENO scheme, comparison of CPU times for computations of discontinuous solution on single-grid and sparse-grid. WENO interpolation for prolongation is used in sparse-grid computations. Final time $T=5/\pi^2$. CFL number is $0.75$. $N_L=3$ in sparse-grid computations. CPU time unit: seconds.}}
\label{ex10_cpu_time}
\end{table}

\section{Conclusions}

In this paper, we develop a third order finite difference  WENO scheme on sparse grids via sparse-grid
combination technique for solving high dimensional hyperbolic problems. Comparable accuracy of the WENO scheme in smooth regions of the solutions
to that of computations on regular single grids is obtained for sparse-grid computations on relatively refined meshes.
A novel WENO prolongation is designed in sparse-grid
combination to achieve the non-oscillatory and high resolution properties of the WENO scheme for resolving shock waves.
With the WENO scheme
on sparse grids, more efficient algorithm than our previous work is obtained for solving the multidimensional hyperbolic equations.
Numerical experiments are performed for the sparse grid WENO method to show significant savings in
computational costs of solving 3D problems by comparisons with single-grid computations.
On relatively coarse meshes, the sparse grid WENO method has larger numerical errors than that by regular single-grid computations.
It will be interesting to improve the accuracy of the sparse grid WENO scheme on coarser meshes, and perform theoretical error analysis for the
scheme.
These will be our future work.

%\newpage

\end{document}